\documentclass[10 pt]{amsart}
\usepackage{amsmath}
\usepackage{amssymb}
\usepackage{amsfonts}
\usepackage{amsthm}
\usepackage{enumerate}
\usepackage[mathscr]{eucal}

\theoremstyle{plain}
\newtheorem{thm}{Theorem}[section]
\newtheorem{lemma}[thm]{Lemma}
\newtheorem{corollary}[thm]{Corollary}
\newtheorem{proposition}[thm]{Proposition}

\numberwithin{equation}{section}

\theoremstyle{plain}
\newtheorem{thmsub}{Theorem}[subsection]
\newtheorem{lemmasub}[thmsub]{Lemma}

\newtheorem{propositionsub}[thmsub]{Proposition}

\numberwithin{equation}{section}

\theoremstyle{remark}

\def\oa{{\overline {\alpha}}}
\def\ob{{\overline {\beta}}}
\def\intslashx{{\frac{1}{|Q_x|}\int_{Q_x}}}
\def\intslash{\rlap{\kern  .32em $\mspace {.5mu}\backslash$ }\int}
\def\qsl{{\rlap{\kern  .32em $\mspace {.5mu}\backslash$ }\int_{Q_x}}}

\def\emph#1{{\it #1 }}
\def\diam{{\text{\rm diam}}}

\def\a{\alpha}

\def\cf{{\it cf}}

\def\dist{{\text{\rm dist}}}

\def\supp{{\text{\rm supp }}}

\def\inn#1#2{\langle#1,#2\rangle}

\def\noi{\noindent}

\def\scr{{s_{\text{cr}}}}

\def\lc{\lesssim}

\def\ga{\gamma}

\def\eps{\varepsilon}

\def\ka{\kappa}
            
\def\la{\lambda}

\def\om{\omega}

\def\fA{{\mathfrak {A}}}
\def\fB{{\mathfrak {B}}}

\def\fJ{{\mathfrak {J}}}

\def\fS{{\mathfrak {S}}}

\def\fm{{\mathfrak {m}}}

\def\bbR{{\mathbb {R}}}

\def\bbZ{{\mathbb {Z}}}

\def\cA{{\mathcal {A}}}
\def\cB{{\mathcal {B}}}
\def\cC{{\mathcal {C}}}

\def\cE{{\mathcal {E}}}
\def\cF{{\mathcal {F}}}

\def\cL{{\mathcal {L}}}
\def\cM{{\mathcal {M}}}

\def\cO{{\mathcal {O}}}

\def\cR{{\mathcal {R}}}
\def\cS{{\mathcal {S}}}
\def\cT{{\mathcal {T}}}

\font \roman = cmr10 at 10 true pt

\def\be#1{\begin{equation}\label{#1}}
\def\bas{\begin{align*}}
\def\eas{\end{align*}}
\def\bi{\begin{itemize}}
\def\ei{\end{itemize}}

\def\dist{{\hbox{\roman dist}}}
\def\supp{{\hbox{\roman supp}}}

\def\eps{\varepsilon}

\def\emph#1{{\it #1}}
\def\textbf#1{{\bf #1}}

\def\intslash{\rlap{\kern  .32em $\mspace {.5mu}\backslash$ }\int}
\def\qsl{{\rlap{\kern  .32em $\mspace {.5mu}\backslash$ }\int_{Q_x}}}
\def\diam{{\text{\rm diam}}}

\begin{document}

\author[M. Pramanik and A. Seeger]
{Malabika Pramanik and Andreas Seeger}
%\title[Averages over curves and associated maximal operators]
%{Averages over curves in $\text{\bf R}^{\text \bf 3}$
%and associated maximal operators}
\title[Averages over curves and associated maximal operators]
%$\text{\bf L}^{\text \bf p}$
{$\mathbf{L}^{\mathbf p}$
Regularity of averages over curves
% in $\text{\bf R}^{\text \bf 3}$
and bounds for  associated maximal operators}

%\date{\today}
%\author{Malabika Pramanik and Andreas Seeger}

\address{M. Pramanik\\
Department of Mathematics\\ University of Rochester \\
Ro\-ches\-ter, NY 14627, USA}

\curraddr{Department of Mathematics, Calfornia Institute of Technology,
Pasadena, CA 91125, USA}

\email{malabika@its.caltech.edu}

\address{A. Seeger   \\
Department of Mathematics\\ University of Wisconsin-Madison\\Madison, WI 53706, USA}
\email{seeger@math.wisc.edu}

\thanks{Version  1-6-05. Supported in part by grants
%DMS 0245408, DMS 0443322 and DMS 0200186
from  the US National Science Foundation.}

\begin{abstract}
We prove that for a finite type curve  in
$\mathbb R^3$ the maximal operator generated by dilations  is bounded
on $L^p$ for sufficiently large $p$. We also show the
endpoint $L^p \rightarrow L^{p}_{1/p}$ regularity result for the
averaging operators for large $p$.
%in this range
%and prove similar results for restricted $X$-ray transforms in $\mathbb
%R^3$ with rigid line complexes.
The proofs make use of a deep result of
Thomas Wolff about decompositions of  cone multipliers.
%\footnote{Include Xray transforms?}
\end{abstract}

\maketitle

\section{Introduction and Statement of Results}\label{intro}

%In this paper we apply Wolff's inequality to derive some new
%consequences concerning maximal functions associated to curves in
%$\mathbb R^3$ and regularity properties of certain averaging
%operators.

 Let $I$ be a compact interval  and consider a smooth curve
$$\gamma\,:\,I \rightarrow \mathbb R^3.$$
%parametrized by
%its arclength $s$, where $I$ is a compact interval. We
%Let $P=\gamma(s_0)$ then
We say that $\gamma$ is of finite type on $I$
if there is a natural number $n$, and $c>0$ so that for all $s\in I$,
and for all $|\xi|=1$,
\begin{equation}
\label{finitetype}
\sum_{j=1}^n\big|\inn{\gamma^{(j)}(s)}{\xi}\big|\ge c.
\end{equation}
For fixed $s$ the smallest
$n$ for which \eqref{finitetype} holds is   the {\it type} of $\gamma$ at $s$.
The type is an upper semicontinuous function, and we refer to the
supremum
of the types over $s\in I$ as the  maximal type of
$\gamma$ on $I$. Let $\chi$ be a smooth function supported in the interior of $I$. We define a measure $\mu_t$ supported on a dilate of the curve by
\begin{equation}  \label{mut}
\langle \mu_t, f \rangle := \int f\left(t {\gamma}(s) \right) \chi(s)\,
ds,
\end{equation}
and set
\begin{equation} \mathcal A_{t}f({x}) := f \ast
  \mu_t({x}).\label{At} \end{equation}
We are aiming to prove  sharp $L^p$ regularity properties of these
integral operators and also $L^p$ boundedness of the maximal operator
given by
\begin{equation}\label{max}
 \mathcal Mf( x):= \sup_{t>0} | \mathcal
  A_{t}f(x)|.
\end{equation}
To the best of our knowledge, $L^p$ boundedness of $\mathcal M$ had not
been previously established for  any $p<\infty$.
Here we prove some positive results for large $p$ and in particular
answer
affirmatively a question on maximal functions associated to helices
which has been around for a while (for example
it  was explicitly
 formulated   in a circulated but  unpublished survey by
Christ \cite{christ-question} from  the late 1980's).

\

Our results rely on a deep inequality of Thomas Wolff for decompositions of
the cone multiplier in $\bbR^3$.
To describe it  consider a
% cone multipliersIn \cite{Wolff1}, T.\ Wolff proved one
%of the most powerful inequalities in Fourier analysis
%which led to essentially sharp $L^p$ local smoothing estimates for the
%wave equation for large $p$.
% Suppose that the Fourier transform of a
distribution $f \in \mathcal S'\left( \mathbb R^3\right)$ whose
Fourier transform is supported
in a
neighborhood of the light cone $\xi_3^2=\xi_1^2+\xi_2^2$ at
level $\xi_3\approx 1$, of width $\delta\ll 1 $.
% approximately $\lambda$, where $\lambda>0$ and $\delta$ is small.
Let $\{\Psi_{\nu}\}$ be a collection of smooth
functions
which are supported in $1 \times
 \delta^{1/2} \times \delta$-plates that
``fit'' the light cone and  satisfy
the  natural size estimates and
differentiability properties;
for a more precise description see
\S\ref{prelim}.  Wolff \cite{Wolff1} proved
that for all sufficiently
large $p$, say $p>p_W$, and all $\epsilon>0$, there exists
$C_{\epsilon,p} > 0$ such that
\begin{equation}
\label{wolffestimate1} \Big\|\sum_{\nu} \widehat \Psi_{\nu} \ast f
\Big\|_p \leq C_{\epsilon,p}
  \delta^{-\frac{1}{2}
+ \frac{2}{p} - \epsilon} \Big( \sum_{\nu} \big\|\widehat{\Psi}_{\nu}
      \ast f \big\|_{p}^{p}\Big)^{\frac{1}{p}}.
\end{equation}
A counterexample in \cite{Wolff1} shows that this inequality
cannot hold for all $\eps>0$ if  $p<6$  and Wolff obtained
\eqref{wolffestimate1} for  $p\ge 74$ (a slightly better range
can be obtained as was observed by
Garrig\'os and one of the authors \cite{GaSe}).

We note that connections between cone multipliers and
the regularity properties of curves
with nonvanishing curvature and torsion have been used
in various previous papers, first implicitly in the paper by
Oberlin \cite{Oberlin87} who proved sharp $L^p\to L^2$ estimates
for, say, convolutions with measures on the helix $(\cos s, \sin s,s)$;
these were extended in \cite{GS} to more general classes of
Fourier integral operators.  Concerning $L^p$ Sobolev estimates the
$L^p\to L^p_{2/3p}$ boundedness
follows by an easy interpolation argument, but improvements of this
estimate are highly nontrivial.
Oberlin, Smith and Sogge
\cite{Oberlin-Smith-Sogge98} used results by
Bourgain \cite{BourgainCone} and Tao and Vargas \cite{TV} on square-functions associated to cone multipliers  to show that
if $2<p<\infty$ then
the averages for the helix map $L^p$ to the Sobolev space $L^p_\alpha$,
for some $\alpha>2/(3p)$.

We emphasize  that sharp regularity results for {\it hypersurfaces} have
 been obtained from
%(somewhat sophisticated)  \
interpolation
arguments, using  results on damped oscillatory integrals
and an improved $L^\infty$ bound near ``flat parts'' of the surface,
see {\it e.g.}
\cite{CoM}, \cite{SoSt}, \cite{Ios} and elsewhere. However this interpolation
 technique
does not apply to averages
over manifolds with very high codimension, in particular not to  curves
in $\bbR^d$, $d\ge 3$.

Our first result on finite type curves in $\bbR^3$ concerns the
averaging operator $\cA\equiv \cA_1$ in \eqref{At}; it
depends on the optimal
exponent $p_W$ in Wolff's inequality \eqref{wolffestimate1}.

\begin{thm} \label{fixedt}
Suppose that  $\gamma\in C^{n+5}(I)$ is of maximal type $n$, and
suppose that
$$\max\{n, \frac{p_W+2}{2}\}<p<\infty.$$ Then $\mathcal A$ maps  $L^{p}$ boundedly to the Sobolev space
$L^{p}_{1/p}$.
\end{thm}
\noindent Thus the sharp $L^p$-Sobolev  regularity properties for
the helix hold for $p>38$, according to Wolff's result. It is known by
an example due to Oberlin and Smith \cite{Oberlin-Smith99} that
the $L^p\to L^p_{1/p}$ regularity result fails if $p<4$. Recall that
 Wolff's inequality \eqref{wolffestimate1} is conjectured for
$p\in (6,\infty)$, and thus establishing this conjecture would by
 Theorem  \ref{fixedt} imply the $L^p\to L^p_{1/p}$ bound for $p>4$.
If the type $n$ is sufficiently large then our result is sharp;
it can be shown by a modification of an  example by Christ \cite{christex}
for plane
curves
that the endpoint
 $L^n\to L^n_{1/n}$ bound fails.
Finally we note that
by a  duality argument and standard facts on Sobolev spaces
 one can also deduce sharp bounds near $p=1$, namely if
$1<p<\min \{n/(n-1),(p_W+2)/p_W\}$
then  $\mathcal A$ maps $L^{p}$
boundedly to $L^{p}_{1/p'}$.

\medskip

\noi{\it Remark.} There are  generalizations of Theorem \ref{fixedt}
which
apply to variable curves; one assumes that the associated
canonical relation in $T^*\bbR^3\times T^*\bbR^3$
projects to each $T^*\bbR^3$  only with fold singularities and that a
 curvature condition  in \cite{GS}
on the fibers of the singular set
is satisfied. We intend to take up these matters
in a forthcoming paper \cite{Pramanik-Seeger}.

\medskip

Our main result on  the maximal operator $\cM$ is

\begin{thm} \label{maximal}
Suppose that $\gamma\in C^{n+5}$ is of maximal type $n$,
then $\mathcal M$ defines a bounded
operator on $L^p$ for $p > \max(n, (p_W+2)/2)$.
\end{thm}

Again the range of $p$'s is only optimal if the maximal
type is sufficiently large ({\it i.e.} $n\ge (p_W+2)/2$).
The following measure-theoretic
consequence (which only uses $L^p$ boundedness for some $p<\infty$)
appears to be new; it follows from Theorem \ref{maximal}
 by  arguments in \cite{Bourgain1991}.

\begin{corollary} Let $\gamma:I\to \bbR^3$  be smooth and  of finite type
and let $A\subset \bbR^3$ be a set of positive measure.
Let $E$ be a subset of $\mathbb R^3$ with the property that for every
$x\in A$ there is a $t(x)>0$ such that $x+t(x){\gamma}(I)$ is
contained in $E$. Then $E$ has positive outer measure.
\end{corollary}

In itself the regularity result of  Theorem \ref{fixedt}
does not imply boundedness of the maximal
operator, but a  local smoothing estimate can be used.
This we only formulate for the  nonvanishing curvature and torsion case.

\begin{thm} \label{localsmoothing}
Suppose that $\gamma\in C^5(I)$ has  nonvanishing curvature and torsion.
Let $p_W<p<\infty$ and $\chi\in C^\infty_0((1,2))$.
%t \in [1,2]$.
Suppose that $\alpha<
4/(3p)$. Then the operator $\fA$ defined by $\fA f({x},t) =
\chi(t) \mathcal A_tf({x})$ maps $L^p(\mathbb R^3)$ boundedly into
$L^p_{\alpha}(\mathbb R^4)$.
\end{thm}
%{\bf{Remark :}}
By interpolation with the standard $L^2 \rightarrow L^{2}_{1/3}$
estimate (obtained from van der Corput's lemma)
one sees that $\fA$ maps $L^p(\mathbb
R^3)$ to $L^{p}_{\beta(p)}( \mathbb R^4)$, for some $\beta(p)
> 1/p$, if $(p_W+2)/2<p<\infty$. By standard arguments the $L^p$
boundedness of the maximal operator $\cM$ follows in this range (provided
that the curve has nonvanishing curvature and torsion).

\medskip

{\it Structure of the paper.}
In \S\ref{prelim} we prove an extension of Wolff's estimate to general cones
which will be crucial for the arguments that follow.
In \S\ref{averages} we prove the sharp $L^p$ Sobolev estimates for large $p$
(Theorem \ref{fixedt}).
In \S\ref{LSM1} we prove a version of the local smoothing estimate
for averaging operators associated to curves in $\bbR^d$ microlocalized to the
nondegenerate region where $\inn{\gamma''(s)}{\xi}\neq 0$.
In \S\ref{LSM2} we use the previous  estimates and
rescaling arguments  to prove
 Theorem \ref{localsmoothing} and in \S\ref{maximalfunctions} we
deduce our results for maximal operators,
including an estimate for a two parameter
family of helices.

%\footnote{Still have to decide whether to include $X$-ray estimate.}

\section{Variations of Wolff's inequality} \label{prelim}

The goal of this section is to prove a variant of
Wolff's estimate \eqref{wolffestimate1}
 where the standard light-cone is replaced by a general cone with one
nonvanishing principal curvature.
Rather than redoing the very complicated proof  of Wolff's
inequality we
shall
use rescaling and induction on scales arguments to
 deduce the general result from the special result, assuming
the validity of \eqref{wolffestimate1} for the light cone, in the range
$p\ge p_W$.

We need to first set up appropriate notation.
Let $I$ be a closed subinterval of $[-1,1]$ and let
\begin{equation*}
\a \mapsto g(\a)=(g_1(\a),g_2(\a))\in \bbR^2, \quad \alpha\in I,
\end{equation*}
define  a $C^3$ curve in the plane and we assume that for positive $b_0$,
$b_1$ and $b_2$
\begin{equation}\label{g-assumption}
\begin{gathered}
\|g\|_{C^3(I)}\le b_0,
\\
|g'(\a)|\ge b_1,
\\
|g_1'(\a)g_2''(\a)-g_2'(\a)g_1''(\a)|\ge b_2.
\end{gathered}
\end{equation}

We consider  multipliers supported near the cone
%in $\bbR^3$
$$\cC_g=\{\xi\in \bbR^3: \xi=\la(g_1(\a), g_2(\a),1), \quad \a \in I, \la>0\}.
$$
For each $\a$ we set
\begin{equation}\label{uialpha}
u_1(\a)= (g(\a),1),\quad u_2(\a)=(g'(\a),0),\quad
u_3(\a)= u_1(\a)\times u_2(\a)
\end{equation}
where $\times$ refers to the usual cross product
so that  a basis of the tangent space of $\cC_g $ at $(g(\a),1)$ is given by
$\{u_1(\a), u_2(\a)\}$.
Let $\la\ge 0$, $\delta>0$ and define the $(\delta,\la)$-plate at $\a$,
$R^\a_{\delta,\la}$,
to be the parallelepiped in $\bbR^3$ given by
the inequalities
\begin{equation}\label{definitionofplates}
\begin{gathered}
\la/2\le |\inn{u_1(\a)}{\xi}|\le 2\la,
%\quad
\\
|\inn{u_2(\a)}{\xi-\xi_3u_1(\a)
}|\le \la\delta^{1/2},
%\quad
\\
|\inn{u_3(\a)}{\xi}|\le \la\delta.
\end{gathered}
\end{equation}
For a constant $A\ge 1$ we define  the {\it $A$-extension} of the plate
$R^\a_{\delta,\la}$  to be the parallelepiped
given by
the inequalities
\begin{equation*}
\begin{gathered}
\la/(2A)\le |\inn{u_1(\a)}{\xi}|\le 2A\la;
\\
|\inn{u_2(\a)}{\xi-\xi_3u_1(\a)}|\le A\la\delta^{1/2};
\\
|\inn{u_3(\a)}{\xi}|\le A\la\delta.
\end{gathered}
\end{equation*}

Note that the $A$-extension of a $(\delta,\la)$-plate has width
$\approx A\la$ in the radial
direction tangent to the cone, width $\approx A\la\delta^{1/2}$ in the tangential direction which is
perpendicular to the radial direction and is supported in a neighborhood of width
$\approx A\la\delta$ of the cone.

A $C^\infty$ function $\phi$ is called an \emph{
admissible bump function associated
to $R^\a_{\delta,\la}$} if $\phi$ is supported in $R^\a_{\delta,\la}$
and if
\begin{multline}\label{estimateforbumps}
\big|
\inn{u_1(\a)}{\nabla}^{n_1}
\inn{u_2(\a)}{\nabla}^{n_2}
\inn{u_3(\a)}{\nabla}^{n_3}
\phi(\xi)\big|\le \la^{-n_1-n_2-n_3} \delta^{-n_2/2}\delta^{-n_3},
\\
\qquad 0\le n_1+n_2+n_3\le 4.
\end{multline}
%\footnote{We may incorporate the shear and
%thus perhaps remove an assumption on $\gam(\a)$.}
A $C^\infty$ function $\phi$ is called an \emph{
admissible bump function associated
to the $A$-extension of  $R^\a_{\delta,\la}$} if $\phi$ is supported in the
$A$-extension but still satisfies the estimates
 \eqref{estimateforbumps}.

Let $\la\ge 1$, $\delta>0$, $\delta^{1/2}\le \theta$, moreover $\sigma\le \delta^{1/2}$.
A finite collection $\cR=\{R_\nu\}_{\nu=1}^N$ is called a
\emph{$(\delta,\la,\theta)$-plate family with separation $\sigma$
associated to $g$}
if (i) each $R_\nu$ is of the form $R^{\a_\nu}_{\delta,\la}$ for some
$\a_\nu\in I$, (ii)  $\nu\neq\nu'$ implies that
$|\alpha_\nu-\alpha_{\nu'}|\ge \sigma$ and, (iii)
$\max_{R_\nu\in \cR}
\{\alpha_\nu\}-\min_{R_\nu\in \cR}\{ \alpha_\nu\}\le \theta.$
Given  $A\ge 1$ we let the {\it $A$-extension of the plate family $\cR$} consist of the $A$-extensions of the plates $R_\nu$.

The main result in Wolff's paper is proved for
the  cone generated by $g(\a)=(\cos 2\pi\a,\sin 2\pi\a)$,
$-1/2\le \alpha\le 1/2$.
Namely if $\cR$ is a
$(\delta,\la,1)$-plate family with separation $\sqrt \delta$ and
for $R\in \cR$, $\phi_R$ is an admissible bump function associated to $R$
then for all $\eps>0$
there is the inequality
\begin{equation}\label{wolffestimate2}
\Big\|\sum_{R\in \cR} \cF^{-1}[\phi_R \widehat f_R]\Big\|_p
\le A(\eps) \delta^{\frac 2p-\frac 12-\eps}
\Big(\sum_{R\in \cR}\|f_R\|_p^p\Big)^{1/p}
\end{equation}
if $p>74$.
This is equivalent with the statement \eqref{wolffestimate1}
in the introduction.
%We also note that scaling  arguments in \cite{TV}, \cite{Wolff1}
%which will be reproduced below
%show
%that if we add the additional restriction that the plates
%lie in a sector of aperture $\theta\ge \delta^{1/2}$ (in other words if
%$$\cR$ is a
%$(\delta,\la,1)$-plate family with separation $\sqrt \delta$) then
%the constant
%$\delta^{\frac 2p-\frac 12-\eps}$ can be improved to
%$(\delta/\theta^2)^{\frac 2p-\frac 12-\eps}$.
%Finally if we relax the $\sqrt \delta$ separation condition
%to a $\sigma$ separation condition for $\sigma<\delta^{1/2} $
%then we get an additional factor  of $\delta^{1/2}/\sigma$;
%this is immediate by the triangle inequality.
Our next proposition
says that this inequality  for the light cone implies an  analogous inequality
for  a general curved cone.

\begin{proposition}\label{generalcones}
Suppose that  $p>2$ and that
 \eqref{wolffestimate2} holds for all
$(\delta,\la,1)$-plate families associated to
 the circle $\{(\cos(2\pi\a),\sin (2\pi\a))\}$.
Then for any  $\delta\le 1$, $\la\ge 1$, $\sigma\le \sqrt\delta$
the following holds true.

Let $\a\mapsto g(\a)$ satisfy
\eqref{g-assumption} and let
$\cR$
%=\{R_\nu\}_{\nu=1}^N$
be  a
$(\delta,\la,\theta)$-plate family with separation $\sigma$, associated to $g$. For each $R$ let $\phi_R$ be an admissible bump function
associated to $R$. Then for $\eps>0$ there is a constant $C(\eps)$ depending
only on $ \eps$ and the constants  $b_0$, $b_1$, $b_2$  in
\eqref{g-assumption} so that for $f_R\in L^p(\bbR^3)$,
% we have the inequality
\begin{equation}\label{wolffestimate3}
\Big\|\sum_{R\in \cR} \cF^{-1}[\phi_R \widehat f_R]\Big\|_p
\le C(\eps) \delta^{1/2}\sigma^{-1}
(\delta\theta^{-2})^{\frac 2p-\frac 12-\eps}
\Big(\sum_{R\in \cR}\|f_R\|_p^p\Big)^{1/p}.
\end{equation}
\end{proposition}

\begin{proof}

 We first remark that
we can immediately reduce to the case $\sigma=\sqrt\delta$, by a
pidgeonhole argument and an application of the triangle inequality.

Secondly, if $\Psi_R$ are bump-functions contained in the
$A$-extensions of the rectangles $R$, but satisfying the
same estimates \eqref{estimateforbumps} relative to the rec\-tangles $R$,
then  an estimate such as
\eqref{wolffestimate3} implies a similar estimate for
the collection of bump functions $\{\Psi_R\}$ where   the constant
$C(\eps)$ is replaced with $C_A C(\eps)$.
This
observation will be used extensively; it is proved by a
pidgeonhole and partition of unity argument.

We now use various scaling arguments based
% of course
on the formula
\begin{equation}\label{scaling}
 \cF^{-1}[m(L\cdot) \widehat f](x)=
\cF^{-1}[ m \widehat{f(L^*\cdot)}] ({L^*}^{-1}x)
\end{equation}
for any real  invertible linear transformation $L$ (with transpose $L^*$).

\medskip

\noi{\it Step 1.} Here we still assume that
$g(\a)=(\cos 2\pi\a,\sin 2\pi\a)$, but wish to show for
$\delta^{1/2}\le \theta\le 1$ the improved estimate
for a
$(\delta,\la,\theta)$-plate family with separation $\delta^{1/2}$.
The plates in $\cR\equiv\{R_\nu\}_{\nu=1}^N$ are of the form
$R_\nu= R^{\a_\nu}_{\delta,\la}$  where $|\a_\nu-\a^0|\le \theta$ for some
$\a^0\in [0,1]$. Let $L_1$ be the rotation by the angle $\a$ in the $(\xi_1,\xi_2)$ plane which leaves the $\xi_3$ axis fixed.
Then the family
the plate family $L_1(R_\nu)$ is still a
$(\delta,\la,\theta)$ plate family  with associated bump functions
%and for suitable $C_0$ the functions
$\phi_{\nu,1}=\phi_{R_\nu}\circ L_{1}^{-1}$.
We note that all rotated plates are contained  in a larger
$(C_1\la, C_1\la \theta, C_1\la\theta^2 )$ rectangle with axes in the direction of
$(1,0,1)$, $(0,1,0)$, $(1,0,-1)$.

We now use a rescaling argument from \cite{TV} and \cite{Wolff1}.
Let $L_2$  be the linear transformation that maps
$(1,0,1)$ to $(1,0,1)$, $(0,1,0)$ to $\theta^{-1}(0,1,0)$ and
$ (1,0,-1)$ to $\theta^{-2}(1,0,-1)$;
%\[(1,0,1) \mapsto (1,0,1), \;
%(0,1,0) \mapsto \frac{1}{\theta}(0,1,0);
% (1,0,-1) \mapsto \frac{1}{\theta_0^2}(1,0,-1), \; \]
it leaves the light cone invariant.
One checks that each parallelepiped $L_2\circ L_1 R_\nu$
is contained in a
$(C_2\la, C_2\la\delta^{1/2}\theta^{-1}, C_2 \la\delta\theta^{-2})$ plate
$\widetilde R_\nu$ and the sets $\widetilde R_\nu$ form a $C_2$ extension
of a $(\delta\theta^{-2}, \la, 1)$ family with separation
$\sigma= C_3^{-1}\delta^{1/2}\theta^{-1}$.
Thus using \eqref{scaling} we may apply the assumed result for $\theta=1$,
and obtain the claimed result for $\theta<1$ (yet for the light cone).

\medskip

\noi{\it Step 2.} We shall now consider tilted cones where $g$ is given by
%where we assume that
\begin{equation}\label{tilt}
g(\a)=(a+\rho \cos\a, b+\rho \sin\a),\qquad |a|+|b|+\rho+1/\rho\le K
\end{equation}

Suppose that  we are given a $(\delta,\la,\theta)$-plate family
$\cR=\{R_\nu\}$  associated to $g$, with
separation $\sqrt\delta$; moreover we are given a family of admissible
bump-functions $\phi_\nu$ associated with the plates $R_\nu$.
Consider the linear transformation $L$ given by
$$L(\xi)=\Xi, \quad\text{ with }
\Xi_1=\frac{\xi_1-a\xi_3}{\rho}, \quad
\Xi_2=\frac{\xi_2-b\xi_3}{\rho}, \quad\Xi_3=\xi_3.
$$
Then the parallelepipeds  $L(R_\nu)$ are contained in
parallelepipeds  $\widetilde R_\nu$ which for a  suitable constant $C_4$
form a  $C_4$-extension
of a
$(\delta,\la,\theta)$-plate family associated to the unit circle;
moreover,
for suitable $C_5$
the functions $C_5^{-1}\phi_\nu\circ L^{-1}$ form an
admissible  collection of bump functions associated to this extension.
Here $C_4$, $C_5$  depend only on
the constant $K$ in \eqref{tilt}.
By scaling we obtain then estimate \eqref{wolffestimate3} for
$g$ as in \eqref{tilt}, with $C(\eps)$ equal to $C(K) A(\eps)$.

\medskip

\noi
{\it Step 3.} We now use an induction on scales argument.
Let $\beta> (1/2-2/p)$ and let
 $W(\beta)$ denote the statement that the inequality
\begin{equation}\label{wolffestimate4}
\Big\|\sum_{R\in \cR} \cF^{-1}[\phi_R \widehat f_R]\Big\|_p
\le B(\beta)
(\delta\theta^{-2})^{-\beta}
\Big(\sum_{R\in \cR}\|f_R\|_p^p\Big)^{1/p}.
\end{equation}
holds for all
$g$ satisfying
\eqref{g-assumption},
all $\la>0$, $\delta\le 1$, $\sqrt{\delta}\le \theta\le 1$, all
$(\delta,\la, \theta)$-plate-families associated to such $g$.

We remark that clearly  $W(\beta)$ holds with
$\beta=1$, with $B(1)$ depending only on the constants in
\eqref{g-assumption}.
%; a straightforward interpolation using the orthogonality this can be
%improved to $\alpha=1-2/p$.
We shall now show that for $\beta>1/2-2/p$
\begin{multline}\label{inductionstep}
 W(\beta)\implies W(\beta')
\quad\text{ with } \\ \beta'=
\tfrac 23\beta+
\tfrac 13(\tfrac 12-\tfrac 2p+\eps),
\quad B(\beta')= C_6 B(\beta) A(\eps),
\end{multline}
where $C_6$  depends only on \eqref{g-assumption}.

In order to show \eqref{inductionstep} we let
$\cR=\{R^{\alpha_\nu}_{\la,\delta}\}$
 be a $(\delta,\la,\theta)$-plate family, with associated family of  bump functions $\{\phi_\nu\}$.
We regroup the indices $\nu$ into families $J_\mu$, so that for
$\nu,\nu'\in J_\mu$ we have $|\alpha_\nu-\a_{\nu'}|\le \delta^{1/3}$ and for
$\nu\in J_\mu$, $\nu'\in J_{\mu+2}$ we have $\alpha_{\nu'}-\alpha_\nu\ge
\delta^{1/3}/2$. For each $\mu$ we pick one $\nu(\mu)\in J_\mu$. Then for
all $\nu\in J_\mu$ the
$R^{\alpha_\nu}_{\la,\delta}$ are contained in the $C_7$-extension $R_\mu'$
of a
$(\delta^{2/3},\la)$-plate  at $\alpha_\nu(\mu)$,
as can be verified by a Taylor expansion. Let $R_\mu''$ be the
$2C_7$ extension of that plate.
We may pick a $C^\infty$-function $\Psi_\mu $ supported in $R_\mu''$ which equals $1$ on $R_\mu'$, so that for suitable $C_8$ depending only on the
constants in
\eqref{g-assumption}, the functions $C_8^{-1}\Psi_\mu$ are admissible bump functions associated to the $R_\mu''$. We then use
assumption $W(\beta)$  to conclude
that
\begin{align}
\Big\|\sum_\nu\cF^{-1}[\phi_\nu\widehat f_\nu]\Big\|_p
&=
\Big\|\sum_\mu\sum_{\nu\in J_\mu} \cF^{-1}[\Psi_\mu\phi_\nu\widehat f_\nu]\Big\|_p
\notag
\\
\label{firstininductionstep}
&\le C_8 B(\beta) (\delta^{2/3}\theta^{-2})^{-\beta}
\Big(\sum_\mu\Big\|\sum_{\nu\in J_\mu} \cF^{-1}[\phi_\nu\widehat f_\nu]\Big\|_p^p\Big)^{1/p}.
\end{align}
We  claim that for each $\mu$,
\begin{equation}\label{secondininductionstep}
\Big\|\sum_{\nu\in J_\mu} \cF^{-1}[\phi_\nu\widehat f_\nu]\Big\|_p
\le C_9 A(\eps) (\delta^{1/3})^{2/p-1/2-\eps}
\Big(\sum_{\nu\in J_\mu}
\|f_\nu\|_p^p\Big)^{1/p}.
\end{equation}
Clearly a combination of
\eqref{firstininductionstep} and
\eqref{secondininductionstep} yields \eqref{inductionstep} with $C_6=C_8C_9$.

 We fix $\alpha_\mu':=\alpha_{\nu(\mu)}$ and observe that
on the interval
$[\alpha_\mu'-\delta^{1/3},\alpha_\mu'+\delta^{1/3}]$
we may approximate
the curve $\alpha\to g(\a)$ by its osculating circle
with accuracy $\le C_{10}\delta$. The circle is   given by
$$
g_\mu(\a)= g(\a_\mu')+ \rho n(\alpha_\mu')+
\rho\big
(\cos (\tfrac{\a - \a_\mu' - \varphi_\mu}{\rho} ),
\sin (\tfrac{\a - \a_\mu' - \varphi_\mu}{\rho} )\big)
$$
where $n(\a)$ is the unit normal vector $(-g_2'(\a),g_1'(\a))/|g'(\a)|$,
 $\rho$ is the reciprocal of the curvature of $g$ at $\alpha_\mu'$
 and $\varphi_\mu$ is the unique value between 0 and $2\pi$
for which
$ |g'(\alpha_\mu')|
\sin(\varphi_\mu/\rho) = g_1'(\a_\mu')$ and
$|g'(\alpha_\mu')|\cos(\varphi_\mu/ \rho) = g_2'(\a_\mu)$.

In view of the good approximation property we see that
for each $\nu\in J_\mu$
the plate $R^{\a_\nu}_{\delta,\la}$ associated to $g$
is contained in the
$C_{11}$-extension $\widetilde R^\nu$ of a plate
$\widetilde R^{\a_\nu}_{\delta,\la}$ associated to
$g_\mu$. Moreover the family $J_\mu$ can be split into no
more than $C_{12}$
subfamilies $J_\mu^i$ where the $\alpha_\nu$ in each subfamily are
$\sqrt{\delta}$-separated. Finally there is $C_{13}$ so that each
 bump function  $\phi_\nu$ is the $C_{13}$-multiple of an admissible bump
function associated to $\widetilde R^\nu$.
Here $C_{11}, C_{12}, C_{13}$ depend only on the constants in
\eqref{g-assumption}. This puts us in the position to apply the result from step 2, with
$\theta=C_{14}\delta^{1/3}$; we observe that the constant $K$ in step 2
controlling in particular  the radius of curvature  depends again
only on the  constants in \eqref{g-assumption}. Thus we can deduce
\eqref{secondininductionstep} and the proof of \eqref{inductionstep}
is complete.

\medskip

\noi{\it Step 4.}  We now iterate \eqref{inductionstep} and
replace $\eps$ by $\eps/2$ to  obtain \eqref{wolffestimate4}
with
\begin{align*}
\beta\equiv \beta_n&=(\tfrac 23)^n+(1-(\tfrac 23)^n)
(\tfrac 12-\tfrac 2p+\tfrac{\eps}2)
\\
B(\beta_n)&= \big(C_{15} A(\tfrac \eps 2)\big)^n, \quad
n=1,2,\dots
\end{align*}
The conclusion of the proposition  follows if we choose
$n>\log(2/\eps)/\log(3/2)$.
\end{proof}

One can use Proposition \ref{generalcones}  and  standard arguments to see
 that results on the circular
cone multiplier in \cite{Wolff1} carry over  to more general cones.
To formulate such a result let $\rho\in C^4(\bbR^2\setminus\{0\})$
be positive away from the origin, and homogeneous of degree $1$.
Consider the Fourier multiplier in $\bbR^3$,
given by $$m_{\la}(\xi',\xi_3)=(1-\rho(\xi'/\xi_3))^\la_+.$$
As in \cite{Wolff1} we obtain

\begin{corollary} Assume that the unit sphere
$\Sigma_\rho=\{\xi'\in \bbR^2:\rho(\xi')=1\}$
has nonvanishing curvature everywhere. Then $m_\lambda$
is a Fourier multiplier of $L^p(\bbR^3)$ if $\la>(1/2-2/p)$,
$p_W\le p<\infty$.
\end{corollary}

\

\emph{Remarks.}

(i) The curvature condition on $\Sigma_\rho$
 in the corollary can be
 relaxed by scaling arguments.

(ii) The methods of Proposition \ref{generalcones} apply
in
higher dimensions as well. In particular they  generalize
Wolff's inequality for decompositions of light cones
 in higher dimensions (\cite{LW}) to more general elliptical cones
generated by convex hypersurfaces with nonvanishing curvature.
In particular, if $\rho$ is a sufficiently smooth
distance function in $\bbR^d$,
$m_\la(\xi',\xi_{d+1})=(1-\rho(\xi'/|\xi_{d+1}|))_+^\la$ and
if the unit sphere associated with $\rho$ is a
convex hypersurface of $\bbR^d$ with nonvanishing Gaussian curvature
 then $m_\lambda$ is a Fourier multiplier
of $L^p(\bbR^{d+1})$, for $\lambda>d|1/2-1/p|-1/2$,
for the range of $p$'s  given in \cite{LW} for the multipliers associated
with the spherical cone.
After   Proposition \ref{generalcones} had been obtained
 \L aba and Pramanik \cite{LP} worked out an alternative
  proof of the higher dimensional variant directly
 based on the methods in \cite{Wolff1},  \cite{LW}. Their
 approach  also applies  to the case of nonelliptical cones
which cannot be obtained by scaling and approximation
 from existing results.

\section{$L^p$ regularity}\label{averages}

We shall first consider a ``nondegenerate case'', namely we assume that
$s\mapsto \gamma(s)\in \bbR^3$, $s\in I\subset [-1,1]$  is of class $C^5$ and
%parametrized by arclength
has {\it nonvanishing curvature and torsion.}
We assume that
\begin{equation} \label{gammaassumption1}
\sum_{i=1}^5|\gamma^{(i)}(s)|\le C_0, \qquad s\in I
\end{equation}
and
\begin{equation} \label{gammaassumption2}
\Big| \det \begin{pmatrix} \gamma'(s)&\gamma''(s)&\gamma'''(s)\end{pmatrix}\Big|\ge c_0, \qquad s\in I
\end{equation}
In this case we show  Theorem \ref{fixedt} under the assumption
that the cutoff function $\chi$ in
\eqref{mut}
is of class $C^4$; then we may without loss of
generalization assume that $\gamma$ is parametrized by arclength
(since reparametrization introduces just  a different $C^4$ cutoff).
In the end of this section we shall describe how to extend the result to the finite type case.

In what follows we shall write
$\cE_1\lc \cE_2$  for two quantities $\cE_1$, $\cE_2$
if
$\cE_1\le C\cE_2$ with a constant $C$ only depending on the constants in
\eqref{gammaassumption1}, \eqref{gammaassumption2}.
We  denote by $T(s), N(s), B(s)$ the Frenet frame of unit tangent, unit normal and unit binormal vector.
We recall the Frenet equations $T'=\ka N$, $N'=-\ka T+\tau B$, $B'=-\tau N$,
with curvature $\ka$ and $\tau$ .
The assumption of nonvanishing curvature and torsion implies that the cone generated by the binormals, $\fB=\{rB(s): r>0, s\in I\}$, has one
nonvanishing principal curvature which is equal
to $r\ka(s)\tau(s)$ at $\xi=rB(s)$.

By localization in $s$ and possible rotation   we may assume
for  the third component of $B(s)$
that
$B_3(s)>1/2$ for all  $s\in I$. If
\begin{equation}\label{gforcone}
g(s)=\big(\tfrac{B_1(s)}{B_3(s)},\tfrac{B_2(s)}{B_3(s)}\big)
\end{equation}
parametrizes the level curve at height $\xi_3=1$ then the curvature property of the cone can be expressed in terms of the curvature of this level curve and
a computation gives
$$\det \begin{pmatrix}
g_1'(s)& g_2'(s)
\\
g_1''(s)&g_2''(s)
\endpmatrix=
\frac{1}{(B_3(s))^3}\det\pmatrix
B_1'(s)&B_2'(s)&B_3'(s)
\\
B_1''(s)&B_2''(s)&B_3''(s)
\\
B_1(s)&B_2(s)&B_3(s)
\end{pmatrix}
=\frac{\kappa(s)\tau(s)}{(B_3(s))^3}.
$$
Thus
%assuming \eqref{gammaassumption1}, \eqref{gammaassumption2}
 the
hypotheses on $g$ in \eqref{g-assumption} are satisfied with constants depending only on
the constants in
\eqref{gammaassumption1}, \eqref{gammaassumption2}.

We shall work with standard Littlewood-Paley cutoffs, and make decompositions of the Fourier multiplier associated to the averages.
Observe first that the contribution of the multiplier near the origin is
irrelevant in view of the compact support of the kernel.
 Thus consider
for $k>0$ the
Fourier multipliers
\begin{equation} \label{dyadicmultipliers}
m_k(\xi)=
\int e^{-it\inn{\gamma(s)}{\xi}} a_k(s,2^{-k} \xi) ds
\end{equation} where we assume that $a_k$
vanishes  outside the annulus
$\{\xi:1/2<|\xi|<2\}$ and satisfies the estimates
\begin{equation} \label{symbolassumptions}
|\partial^j_s \partial^{\alpha}_\xi a(s,\xi)|\le C_2,\quad
|\alpha|\le 2, 0\le j\le 3;
\end{equation}
here of course $|\alpha|=|\alpha_1|+|\alpha_2|+|\alpha_3|$.
%\footnote{we may need more.}
Thus the multiplier $m_k$   is a symbol of order $0$, with perhaps limited order of differentiability,  localized to the
annulus $\{|\xi|\approx 2^k\}$). We note that by the standard Bernstein theorem
(which says  that $L^2_{\alpha}\subset \cF[L^1]$ for $\alpha>3/2$)
the multipliers  $a_k(s,\cdot )$ and their $s$-derivatives up to order three
are  Fourier multipliers of $L^p(\bbR^3)$, uniformly in $s$, $k$.

We have to establish that for the desired range of $p$'s the sum
$\sum_{k>0} 2^{k/p} m_k$ is a Fourier multiplier of $L^p$.
We may assume that the symbols $a_k$ are supported near from the cone
generated by the binormal vectors $B(s)$. More precisely if
$\theta(\xi)$ is smooth away from the origin and  homogeneous  of
degree $0$ and if
$\theta$  has the property
that $$|\inn{\gamma''(s)}{\tfrac{\xi}{|\xi|}}|\ge c>0, \qquad \xi\in\supp(\theta)\cap
\supp (a_k),
$$ then
$\|\theta m_k\|_\infty=O(2^{-k/2})$ by van der Corput's Lemma, and by
the almost disjointness of the supports
we also have that
$\|\theta\sum_{k>0} 2^{k/2} m_k\|_\infty=O(1
)$ by van der Corput's Lemma.
 Moreover by standard singular integral theory
 the operator with Fourier multiplier
$\theta \sum_{k>0}
m_k$ maps $L^\infty$ to $BMO$ and consequently, by
analytic interpolation
$\theta \sum_{k>0} 2^{k/p}
m_k$ is a Fourier multiplier of $L^p$
 provided $2\le p<\infty$.

Thus by a partition of unity  it suffices to  understand the localization
of the multiplier
$\sum_{k>0} 2^{k/p}m_k$ to a narrow  (tubular)
neighborhood of the binormal cone $\fB=\{rB(s): r>0, s\in I\}$,
and therefore in what follows we may and shall
assume that
 $\xi$ in
 the support of
$a_k(s,\cdot) $ can be expressed  as
$$
\xi= r B(\sigma) +uT(\sigma)=:\Xi(r,u,\sigma),
$$
with inverse  function $\xi\mapsto (r(\xi), u(\xi), \sigma(\xi))$.

\medskip

\noi{\bf Decomposition of the dyadic multipliers.}
We shall now concentrate on the multipliers $m_k$  in
\eqref{dyadicmultipliers}, and prove the bound
$\|m_k\|_{M_p}\lc 2^{-k/p}$, for $p>(p_W+2)/2$.
Here $M_p$ is the usual Fourier multiplier space.

We first decompose further  our symbols $a_k$.
Let $\eta_0\in C^\infty_0(\bbR)$ be an even function supported in
$[-1,1]$ and be equal to $1$ on $[-1/2,1/2]$.
Let $\eta_1=\eta_0(4^{-1}\cdot)-\eta_0$. Let
 $A_0\gg2 \max\{1, 1/\tau(s): s\in I\}$
and set
\begin{equation}\label{definitionofwidetildeak}
\widetilde a_k(s,\xi)=a_k(s,\xi)
\eta_0(2^{2[k/3]}(|u(\xi)|+ (s-\sigma(\xi))^2)),
\end{equation}
and for integers $l<k/3$
\begin{equation}\label{definitionofaklbkl}
\begin{aligned}
a_{k,l}(s,\xi)&= a_k(s,\xi)\,
\eta_1(2^{2l}(|u(\xi)|+ (s-\sigma(\xi))^2))\,
\eta_0(\tfrac{(s-\sigma(\xi))^2}{A_0u(\xi)})
 \\
b_{k,l}(s,\xi)&= a_k(s,\xi)\,
\eta_1(2^{2l}(|u(\xi)|+ (s-\sigma(\xi))^2))\,
\big(1-\eta_0(\tfrac{(s-\sigma(\xi))^2}{A_0 u(\xi)})\big).
\end{aligned}
\end{equation}
Thus $a_{k,l}(s,\cdot)$  is supported where
$\dist(\xi, \fB)\approx 2^{-2l}$
and $|s-\sigma(\xi)|\lc 2^{-l}$, and
$\widetilde a_k(s,\cdot)$
is supported in a $C2^{-2k/3} $
neighborhood of the binormal cone with $|s-\sigma(\xi)|\lc 2^{-k/3}$.
 The symbol $b_{k,l}(s,\cdot)$  is  supported in a
$C2^{-2l}$ neighborhood of the binormal cone  but now
 $|s-\sigma(\xi)|\approx 2^{-l}$.

We note that in view of the preliminary localizations the  symbols
$a_{k,l}$, $b_{k,l}$ vanish for $l\le C$, moreover
$$a_k(s,\xi)=\widetilde a_k(s,\xi)+ \sum_{l\le k/3} a_{k,l}(s,\xi)
+ \sum_{l\le k/3} b_{k,l}(s,\xi).
$$
Set
\begin{equation}\label{mka}
m_k[a](\xi)= \int a(s,2^{k}\xi) e^{-i \inn{\gamma(s)}{\xi}} ds.
\end{equation}
We shall show
\begin{proposition}\label{besov}
For $p_W<p<\infty$,
\begin{align}
\|m_k[\widetilde a_{k}]\|_{M^p} &\le  C_\eps 2^{-4k/3p+k\eps },
\label{tildeak}
\\
\|m_k[a_{k,l}]\|_{M^p} &\le  C_\eps 2^{-k/p}2^{-l/p+l\eps},
\label{akl}
\\
\|m_k[b_{k,l}]\|_{M^p} &\le  C_\eps 2^{-2k/p}2^{2l/p+l\eps},
\label{bkl}
%\|m_k[\widetilde b_{k}]\|_{M^p} &\le  C_\eps 2^{-4k/3p+k\eps }.
%\label{tildebk}
\end{align}
The constants  depend only on $\eps$,
\eqref{gammaassumption1}, \eqref{gammaassumption2} and
\eqref{symbolassumptions}.
\end{proposition}
We shall give the proof of \eqref{akl} and \eqref{bkl}, and the proof of
\eqref{tildeak} is
%\eqref{tildebk}
 are  analogous with mainly notational changes.

For the proofs of
\eqref{akl} and \eqref{bkl}
 we need to further split the symbols $a_{k,l}$, $b_{k,l}$
 by making an equally spaced
decomposition into pieces supported on $2^{-l}$ intervals.
Let $\zeta\in C^\infty_0$ be supported in $(-1,1)$ so that
$\sum_{\nu\in \bbZ}\zeta(\cdot-\nu)\equiv 1$.
We set
$$a_{k,l,\nu}(s,\xi)=\zeta(2^ls-\nu)a_{k,l}(s,\xi)$$
and similarly $b_{k,l,\nu}(s,\xi)=\zeta(2^l s-\nu)b_{k,l}(s,\xi)$; moreover define $\widetilde a_{k,\nu}(s,\xi)=\zeta(2^{k/3}s-\nu) \widetilde a_k(s,\xi)$.

In order to apply Wolff's estimate in the form of Proposition
\ref{generalcones}  we need

\begin{lemma} \label{platesandestimates}
Let $s_\nu= 2^{-l}\nu$ and let $(T_\nu, N_\nu, B_\nu)=
(T(s_\nu),N(s_\nu),B(s_\nu))$.
Suppose that
 $|s-s_\nu|\le 2^{2-l}$.
Then the following holds true:

(i)
The multipliers $a_{k,l,\nu}(s,\cdot)$,
$b_{k,l,\nu}(s,\cdot)$ are supported in
\begin{equation}\label{platedef}
\{\xi:
|\inn{\xi}{T_\nu}|\le C 2^{-2l},
|\inn{\xi}{N_\nu}|\le C 2^{-l},
C^{-1}\le |\inn{\xi}{B_\nu}|\le C\}
\end{equation}
where $C$ only depends on the constants in \eqref{gammaassumption1}.

(ii) For $j=0,1,2$, and $h_\nu=a_{k,l,\nu}(s,\cdot)$ or
$b_{k,l,\nu}(s,\cdot)$
 \begin{align}\label{platediff}
&\big|\big(\inn{T_\nu}{\nabla}\big)^{j}  h_\nu\big|
\le C' 2^{2lj},
\\
&\big|\big(\inn{N_\nu}{\nabla}\big)^{j}  h_\nu\big|
\le C' 2^{lj},
\\
&\big|\big(\inn{B_\nu}{\nabla}\big)^{j}  h_\nu\big|
\le C'.
\end{align}

(iii) The statements analogous to (i), (ii)  hold true
 for $\widetilde a_{k,\nu}(s,\cdot)$,
$\widetilde b_{k,\nu}(s,\cdot)$, with $l$ replaced by $[k/3]+1$.

(iv) If $h_\nu$ is any of the multipliers
 $a_{k,l,\nu}(s,\cdot)$,
$b_{k,l,\nu}(s,\cdot)$,
then the statements analogous to (i)-(iii)
hold for the multiplier
$h_\nu(\xi) 2^l\inn{\gamma''(s)}{\xi}$.
Similarly, if $\widetilde h_\nu$ denotes any of
$\widetilde a_{k,\nu}(s,\cdot)$ or
$\widetilde b_{k,\nu}(s,\cdot)$
then $\widetilde h_\nu$ can be replaced by
 $\widetilde h_\nu 2^{l}\inn{\gamma''(s)}{\xi}$.
\end{lemma}

\begin{proof}
To see the containment of
$\supp \,a_{k,l,\nu}(s,\cdot)$ in the set
\eqref{platedef} we assume that $\xi=\Xi(r,u,\sigma)$ and expand
$B(\sigma)$, $T(\sigma)$ about $\sigma=s_\nu$. Using the
 Frenet formulas for $\xi$ in the support of
$a_{k,l,\nu}(s,\cdot)$ we obtain
$$\inn{\xi}{T_\nu}=
\inn{rB(\sigma)+uT(\sigma)}{T_\nu}=
O((\sigma-s_\nu)^2)+O(u)= O(2^{-2l})
$$
and similarly $\inn{\xi}{N_\nu}=O(2^{-l})$.

To show
 \eqref{platediff} we use the
formulas
\begin{equation*}
\nabla r=B, \quad\nabla u= T, \quad
\nabla\sigma= \frac 1{u\kappa -r\tau} N;
\end{equation*}
here of course  $r=r(\xi)$, $B=B(\sigma(\xi)) $, etc. Moreover
\begin{align*}
\inn{e}{\nabla}^2 r&=\frac{-\tau}{-r\tau+u\kappa} \inn{e}{N}^2,
\\
\inn{e}{\nabla}^2 u&=\frac{1}{-r\tau+u\kappa} \inn{e}{N}\inn{e}{T},
\\
\inn{e}{\nabla}^2 \sigma&=\frac{1}{-r\tau+u\kappa} \inn{e}{N}
\inn{e}{\nabla(\tfrac{1}{-r\tau+u\kappa})}.
\end{align*}
{}From these formulas and the
chain rule the verification of the
asserted differentiability properties is straightforward;
we use  also that $T_\nu-T(\sigma(\xi))=O(2^{-l})$ and similar
statements for $N_\nu$ and $B_\nu$.
\end{proof}

We shall need bounds for the $L^1$ and $L^2$ operator norms of the operators defined by
\begin{equation}\label{defAk}
\widehat {\cA^{k,l,\nu} f}(\xi)= m_k[a_{k,l,\nu}]\widehat f(\xi),
\quad \widehat {\widetilde \cA^{k,\nu} f}(\xi)= m_k[\widetilde a_{k,\nu}]
\widehat f(\xi),
\end{equation}
and
\begin{equation}
\label{defAb}
\widehat {\cB^{k,l,\nu} f}(\xi)= m_k[b_{k,l,\nu}]\widehat f(\xi).
\end{equation}
We remark that part (iii) of the following lemma
(and also part (iv) of Lemma \ref{platesandestimates} above)
is not needed in this section but will be needed in a proof of
Theorem \ref{Kakeya}.

\begin{lemma}\label{L1L2bounds}
(i)
\begin{align}\label{AklnL2}
\|\cA^{k,l,\nu}\|_{L^2\to L^2} &\le C 2^{(l-k)/2},\quad      l\le k/3,
\\
\label{tAknL2}
\|\widetilde \cA^{k,\nu}\|_{L^2\to L^2}
&\le C 2^{-k/3},
\\
\label{BklnL2}
\|\cB^{k,l,\nu}\|_{L^2\to L^2} &\le C 2^{2l-k},\quad      l\le k/3.
%\\\label{tBknL2}
%\|\widetilde \cB^{k,\nu}\|_{L^2\to L^2}
%&\le C 2^{-k/3}.
\end{align}

(ii)
\begin{align}\label{knlLinfty}
\|\cA^{k,l,\nu}\|_{L^\infty\to L^\infty}
+\|\cB^{k,l,\nu}\|_{L^\infty\to L^\infty}
&\le C 2^{-l},\quad      l\le k/3,
\\
\label{tknLinfty}
\|\widetilde \cA^{k,\nu}\|_{L^\infty\to L^\infty}
%+ \|\widetilde \cB^{k,\nu}\|_{L^\infty\to L^\infty}
&\le C 2^{-k/3}.
\end{align}

(iii) Assume now that the number of sign changes of the function
$s\mapsto \inn{\gamma'''(s)}{\xi} $ is bounded independent of $\xi$.
Then the estimates in (i), (ii)  continue to hold true if we replace
in the above definitions any of the symbols
$h_\nu=
a_{k,l,\nu}(s,\cdot)$ or $b_{k,l,\nu}(s,\cdot)$
with $h_\nu 2^{l}\inn{\gamma''(s)}{\xi}$,
or  if we replace
$\widetilde h_\nu=\widetilde  a_{k,\nu}(s,\cdot)$
%or $\widetilde  b_{k,\nu}(s,\cdot)$
with $\widetilde h_\nu 2^{k/3}\inn{\gamma''(s)}{\xi}$.

\end{lemma}

\begin{proof}
The $L^2$ estimates \eqref{tAknL2}
%and \eqref{tBknL2}
are immediate from van der Corput's lemma with third derivatives;
 we use that $\inn{\gamma'''(s)}{2^k \xi}\approx 2^k$ for small
$u(\xi)$.
We use van der Corput's estimate for
\eqref{AklnL2}
as well and observe that for $\xi\in \supp \, a_{k,l,\nu}$
we have that $\inn{\gamma''(s)}{\xi}=(s-\sigma(\xi))
\inn{\gamma'''(s)}{\xi}+ O(2^{-2l})$ so that
$\inn{\gamma''(s)}{2^k\xi}\approx 2^{k-l}$ if
$|s-\sigma(\xi)|\ge c_0 2^{-l}$.
If $c_0$ is  sufficiently small then we also have for
$|s-\sigma(\xi)|\le c_0 2^{-l}$ that
$\inn{\gamma'(s)}{\xi}=
\inn{\gamma'(\sigma(\xi))}{\xi}
+ O(c_0^2 2^{-2l})$ and since
$|\inn{\gamma'(\sigma(\xi))}{\xi}|=|u(\xi)|$ we get
$\inn{\gamma'(s)}{2^k\xi}\approx
2^{k-2l}$ if
$|s-\sigma(\xi)|\le c_0 2^{-l}$.
Thus van der Corput's
lemma with one or two derivatives yields the bound
$$\big\|m_k[a_{k,l,\nu}]\big\|_\infty
\le C(2^{(l-k)/2}+ 2^{2l-k})\le C' 2^{(l-k)/2},
$$ since $l\le k/3$.

A similar argument goes through for
$m_k[b_{k,l,\nu}]$. Now  $|u(\xi)|\ll |s-\sigma(\xi)|\approx 2^{-l}$ so that
in the support of $b_{k,l,\nu}$ there is the
lower bound $|\inn{\gamma'(s)}{2^k\xi}|\ge c 2^{k-2l}$
and van der Corput's lemma with one derivative yields
$$\|m_k[b_{k,l,\nu}]\|_\infty\lc 2^{2l-k}$$
and thus the asserted $L^2$ bound \eqref{BklnL2} .

We now turn to the $L^\infty$ bounds.  Consider first
the multiplier
$a_{k,l,\nu}$. Let $L_\nu$ be the rotation that maps the
coordinate vector $e_1$ to
$T_\nu$, $e_2$ to $N_\nu$ and $e_3$ to $B_\nu$. Let $\delta_l$ denote the nonisotropic dilation defined by $\delta_l(\xi)=(2^{-2l}\xi_1,2^{-l}\xi_2, \xi_3)$.
By scaling we see from
\eqref{platedef} and \eqref{platediff} that
$a_{k,l,\nu}(L_\nu\delta_l \cdot)$
 is supported on a ball of radius $C$ and that the  directional
derivatives up to order  $2$ in the $e_1,e_2,e_3$ directions
are bounded, uniformly in $k,l,\nu,s$. Thus we may apply
Bernstein's theorem (alluded to above after formula \eqref{symbolassumptions})
%which is the embedding of the Sobolev-space
%$L^2_\alpha(\bbR^3)$, $\alpha>3/2$,
% in $\cF^{-1}[L^1]$
and we see that the $L^1$ norms of the
 functions
$\cF^{-1}[a_{k,l,\nu}(s,L_\nu\delta_l \cdot)]$ are uniformly bounded.
By scaling and translation we also see that the $L^1$ norms of the
functions $\cF^{-1}[a_{k,l,\nu}(s,2^k\cdot)e^{i\inn{\gamma(s)}{\cdot}}]$
are  uniformly bounded  and thus
$$
\big\| \cF^{-1}\{ m_k[a_{k,l,\nu}]\}\big\|_1\le
\int_{|s-s_\nu|\le 2^{2-l}}
\big\|\cF^{-1}[a_{k,l,\nu}(s,2^k\cdot)e^{i\inn{\gamma(s)}{\cdot}}]\big\|_1 ds \le C 2^{-l}.
$$
This implies
the claimed $L^\infty$ bound for
$\cA^{k,l,\nu}$. The other estimates in (ii) are obtained in the same way.

Finally we examine the statement in (iii).
We note for the $L^2$ bounds that
$\inn{\gamma''(s)}{\xi}=O(2^{-l})$ in the support of $a_{k,l,\nu}$;
moreover  for $|\xi|\approx 1$  the integral
$\int|\inn{\gamma'''(s)}{\xi}|ds $ (over the support of the relevant
cutoff function) is also $O(2^{-l})$, by an application of the fundamental theorem of calculus to a bounded number of intervals on which
$\inn{\gamma'''(s)}{\xi}$ has constant sign. This estimate  is
needed for the application of van der Corput's lemma as before where we
now gain a factor of $2^{-l}$. A quick examination of the argument in
Lemma
\ref{platesandestimates} gives the claimed $L^\infty$ bounds
for this case.
\end{proof}

\noi\emph{Proof of Proposition \ref{besov}.}
We prove \eqref{akl}.
%Using the notation in \eqref{defAk} we have
Observe
$m_k[a_{k,l}]=\sum_\nu m_k[a_{k,l,\nu}]$ where
the multipliers
$m_k[a_{k,l,\nu}]$ are supported in $C$-extensions of
$(2^k, 2^{-2l})$ plates associated to the cone generated by
$g(s)$ as in \eqref{gforcone}. This family of plates is a union of a
 bounded number of $c2^{-l}$ separated plate families.
Consequently we can apply Wolff's estimate in the form of Proposition
\ref{generalcones} and we get for $p>p_W$
\begin{equation}\label{Step1}
\Big\|\sum_\nu \cA^{k,l,\nu} f\Big\|_p \le C_\eps 2^{2l(\frac 12-\frac 2p+\eps)}
\Big(\sum_\nu \big\|
\cA^{k,l,\nu} f\big\|_p^p\Big)^{1/p}
\end{equation}
Next we claim that for $2\le p\le \infty$
\begin{equation}\label{Step2}
\Big(\sum_\nu \big\|
\cA^{k,l,\nu}
f\big\|_p^p\Big)^{1/p}
 \le C 2^{-l(1-3/p)}2^{-k/p}\|f\|_p
\end{equation}
where for $p=\infty$ we read the left hand side as an
$\ell^\infty(L^\infty)$ norm. The case for $p=\infty$ follows from
\eqref{knlLinfty} and the case $p=2$ follows from
\eqref{AklnL2} if we also use the finite overlap of the supports of the
multipliers $m_k[a_{k,l,\nu}]$. The case for $2<p<\infty$ follows by interpolation. Now the desired bound
\eqref{akl} follows from  \eqref{Step1} and \eqref{Step2}.

The estimate
\eqref{Step1} holds still true if we replace
$\cA^{k,l,\nu}$ by
$\cB^{k,l,\nu}$. Moreover the argument leading to \eqref{Step2} equally applies, except that we now have a better $L^2$ bound
$O(2^{2l-k})$ and consequently the $\ell^p(L^p)$ bound improves to
\begin{equation*}
\Big(\sum_\nu \big\|
\cB^{k,l,\nu}
f\big\|_p^p\Big)^{1/p}
 \le C 2^{-l(1-6/p)}2^{-2k/p}\|f\|_p.
\end{equation*}
This yields
\eqref{bkl} and the proof of the bound \eqref{tildeak} is analogous.\qed

\medskip

By a further interpolation we also obtain

\begin{corollary}\label{besovextended} For $p> (p_W+2)/2$ there is $\eps_0=\eps_0(p)>0$ so that
\begin{align}
\|m_k[a_{k,l}]\|_{M^p}
&\le  C_p 2^{-k/p}2^{-\eps_0 l/p},
\label{aklcor}
\\
\|m_k[\widetilde a_{k}]\|_{M^p} &\le  C_p 2^{-k(1+\eps_0)/p }.
\label{tildeakcor}
\end{align}
Moreover
\begin{align}
\sum_{k\ge 3l}2^{k/p}\|m_k[b_{k,l}]\|_{M^p} &\le  C_p,
\label{bklcor}
\\
\sum_k 2^{k/p} \|m_k[\widetilde b_{k}]\|_{M^p} &\le  C_p.
\label{tildebkcor}
\end{align}
\end{corollary}

\begin{proof} By the almost disjointness of our plate families and
the $L^2$ bounds in Lemma \ref{L1L2bounds}
 we see that
\begin{equation} \label{L-2}
\|m_k[a_{k,l}]\|_{\infty}\le C 2^{(l-k)/2},
\qquad \|m_k[\widetilde a_{k}]\|_{\infty}=O(2^{-k/3})
\end{equation}
and, similarly,
$\|m_k[b_{k,l}]\|_{\infty}=O(2^{(l-k)/2}$ and
$\|m_k[\widetilde b_{k}]\|_{\infty}=O(2^{-k/3})$.
Interpolating the resulting $L^2$ estimates with the $L^{p}$  bounds
of Proposition \ref{besov} yields the assertion.
\end{proof}

\medskip

\noi{\bf Sobolev estimates.} In order to prove Theorem \ref{fixedt}
we will still have to put the estimates
\eqref{aklcor} for different $k$   together.
The  desired estimates for the corresponding expressions
involving
$m[\widetilde a_k]$,
$ m[b_{k,l}]$ and $m[\widetilde b_k]$   follow of course from Corollary
\ref{besovextended}.
To finish the proof of Theorem \ref{fixedt} for the case of nonvanishing curvature and torsion it suffices to show that
\begin{equation*}
\Big\|\sum_{k\ge 3l} 2^{k/p} m[a_{k,l}]\Big\|_{M^p} \le C_p 2^{-l\eps_1(p)}
\end{equation*} with  $\eps_1(p)>0$ if $p>(p_W+2)/2$.
In what follows we define the operator
$\cA^{k,l} $ by
$\widehat {\cA^{k,l} f}
= m[a_{k,l}]\widehat f$.

By Littlewood-Paley theory it is sufficient to prove the
vector-valued inequality
\begin{equation}
\Big\|\Big(\sum_{k\, :\, k \geq 3l}|2^{k/p} \mathcal A^{k,l} f_k|^2\Big)^{1/2}\Big\|_p
\lesssim 2^{-\epsilon_1(p) l}
\Big\|\Big(\sum_{k>0}| f_k|^2\Big)^{1/2}\Big\|_p.
\label{maincz}
\end{equation}
where $\epsilon(p)>0$, $p>(p_W + 2)/2$.

To verify \eqref{maincz} we follow closely an argument in
\cite{SeTAMS} and use  a vector-valued version of
the Fefferman-Stein inequality for the  $\#$-function and linearization.
The result in \cite{SeTAMS} does not apply but the method does if we
replace certain estimates for singular integrals
by $L^\infty\to BMO$ estimates for averaging operators
({\it cf.} the bound for
\eqref{estimateforII}
below).

Let us consider a family of cubes $Q_x $ with $x\in Q_x $
so that the corners of $Q_x $ are measurable functions, and suppose that
\begin{equation}
\label{gk}
\sup_{x,y} (\sum_k |g_k(x,y)|^2)^{1/2}\le 1.
\end{equation}
%In the following argument, the slashed integral
%$\intslash_Q f$ represents the average $|Q|^{-1}\int_Q f$.
We define, for $k\ge 3l$,
the linearized operator
\begin{equation}
T^z_{l,k}f(x)=
 2^{k(1-z)/2}
%\intslash_{Q_x}
\frac{1}{|Q_x|}\int
\Big[
\mathcal A^{k,l} f(y)-
%\intslash_{Q_x}
 \int\mathcal A^{k,l} f(u)\frac{ du}{|Q_x|}\Big]  g_k(x,y) dy,
% \qquad k\ge 3l
\end{equation}
and also
%%$T^z_{l,k}f=0$ if k<3l$; moreover
 \[
T_l^z F(x) =\sum_{k\ge 3l} T^z_{l,k} f_k(x)\quad \text{ where } F = \{f_k \} \in L^p(\ell^2).
\]
The exponent $p(z)$ is given by $1/p(z)=(1-\text{Re}(z))/2$.
We shall have to prove that for $p(z)>(p_W+2)/2$ the operator
 $T_l^z $ maps $L^{p(z)}(\ell^2)$ to  $L^{p(z)}$
 with operator norm
$O(2^{-l\epsilon})$,
independent of  the choice of $Q(\cdot)$ and $g_k(\cdot,\cdot)$.

Split $T^z_{l} F= I^z_{l} F +II^z_{l}  F+ III^z_{l} F$, where
{\allowdisplaybreaks
\begin{align}
I_l^z F(x) &=
\sum_{\substack{
k>0\\ 2^{-10 l}\le 2^k \diam Q_x \le 2^{10 l}
}}
T_{l,k}^z f_k(x),
\\
II_l^z F(x)&=
\sum_{\substack{
k>0\\  2^k \diam Q_x \ge 2^{10 l}
}}
T_{l,k}^z f_k(x),
\\
III_l^z F(x)&=
\sum_{\substack{
k>0\\  2^k \diam Q_x \le 2^{-10 l}
}}
T_{l,k}^z f_k(x).
\end{align}
}
The main term is $I_l^zF(x)$ which is bounded by

\begin{multline*}
%\intslashx
\frac{1}{|Q_x|}
\int_{Q_x} \Big(
\sum_{\substack{
k>0\\ 2^{-10 l}\le 2^k \diam Q_x \le 2^{10 l}
}}2^{2k/p}
\big|
\mathcal A^{k,l}_t f_k(y)-
\intslashx
\mathcal A^{k,l}_t f_k(u)\frac{ du}{|Q_x|}
\big|^2\Big)^{1/2}dy
\\
%& \hskip2in \text{(by Cauchy-Schwarz inequality and (\ref{gk}))} \\
%&\lesssim (1+l)^{1/2-1/p}
% \intslash_{Q_x}
%\Big(
%\sum_{k>0} 2^k \big|
%\mathcal A^{k,l}_t f_k(y)-
%\intslash_{Q_x} \mathcal A^{k,l}_t f_k(u) du\big|^p\Big)^{1/p} dy
%\\
%& \hskip3in \text{(by H\"older's inequality)} \\
\lesssim (1+l)^{1/2-1/p}\Big(
\sum_{k>0} \big[ 2^{k/p}
M_{HL}(
\mathcal A^{k,l} f_k)\big]^p\Big)^{1/p}(x),
\end{multline*}
%}
by H\"older's inequality.
Here $p = p(z)$ and $M_{HL}$ denotes the Hardy-Little\-wood maximal function. Now for $p=p(z)>(p_W+2)/2$,
\begin{align*}
\|I_l^z F\|_{p}
&\lesssim (1+l)^{1/2-1/p} \sup_k 2^{k/p}\|\mathcal A^{k,l}\|_{L^{p}\to L^{p}}
\|F\|_{L^p(\ell^p)},
%&\lesssim  (1+l)^{1/2-1/p}
%\Big(\sum_{k>0}
%\big\| 2^{k/p} \mathcal A^{k,l}_t f_k\|_{p}^{p}\Big)^{1/p}
%\\
%&\lesssim (1+l)^{1/2-1/p} \left(\sup_k 2^{k/p}\|\mathcal
%  R^{k,l}\|_{L^{p}\to L^{p}} \right)
%\Big(\sum_{k>0}
%\big\|f_k\|_{p}^{p}\Big)^{1/p}
\\
&\le C_p (1+l)^{1/2-1/p} 2^{-l\eps_0(p)}\|F\|_{L^p(\ell^2)},
\end{align*}
by Corollary \ref{besovextended}.

For the operators $II_l^z$ and $III_l^z$ we prove $L^2(\ell^2)\to L^2$
boundedness
for $z=i\tau$  and
$L^\infty (\ell^2)\to L^\infty $ boundedness  for $z=1+i\tau$, with
bounds uniform in $\tau$. The $L^p(\ell^p)\to L^p$ estimate for
 $(1-\text{Re} (z))/2=1/p$ then follows by analytic interpolation.

Using orthogonality arguments we obtain
\begin{equation}
\|II^{i\tau}_l F\|_2+\|III^{i\tau}_l F\|_2\le C_{\epsilon}2^{l(\frac{1}{2} + \epsilon)}
\|F\|_{L^2(\ell^2)} \label{II-two}
\end{equation}
for arbitrary $\epsilon >0$. To see this, let us consider  $II^{i \tau}_l$.
By \eqref{gk}
\begin{align*}
&\left| II^{i \tau}_l F (x) \right| \\&\leq
\sum_{\begin{subarray}{c}
   2^k \diam(Q_x) \geq 2^{10l} \end{subarray}}
\frac{2^{k/2}}{|Q_x|}\int_{Q_x}
%%\intslashx
\big| \mathcal A^{k,l} f_k(y) - \intslashx\mathcal
  A^{k,l} f_k(u) \, du \big| |g_k(x,y) | \, dy \\
 &\leq \Big(\sum_{k>0} \big|2^{\frac{k}{2}} \big(M_{HL} (\mathcal A^{k,l}
  f_k )\big)(x)\big|^2\Big)^{1/2}.
\end{align*}

Therefore,
\[\| II^{i \tau}_l F \|_{2} \leq \Big(\sum_{k}\| 2^{\frac{k}{2}}
    \mathcal A^{k,l}f_k\|_2^2\Big)^{1/2}
\leq C 2^{l/2}
\|F\|_{L^2(\ell^2)}, \]
where we have used (\ref{L-2}).
This proves \eqref{II-two} for
$II^{i\tau}_l$ and the argument for $III^{i\tau}_l$ is exactly analogous.

For the $L^\infty$ bounds let us consider $II^{1+i\tau}_l F(x)$ for fixed $x$ and $\text{Re} (z)=1$. We note that
\begin{equation}\label{estimateforII}
II^{1+i\tau}_l F(x)\le
2 % \intslash_{Q_x}
 \frac{1}{|Q_x|}\int_{Q_x}
\Big(\sum_{2^k\diam(Q_x )> 2^{10 l}}
\big|
\mathcal A^{k,l} f_k(y)\big |^2\Big)^{1/2}dy.
\end{equation}
Let
\[
\mathcal U(x)=\{ y: |x-y+\gamma(s)|\le 10\,\diam (Q_x ), \text{ for
  some } s\in \text{supp}(\chi). \}
\]
Then \[|\mathcal U(x)|\lesssim \big(\diam (Q_x) \big)^2.\]
We estimate
\[II^{1+i\tau}_l F(x)\le
2\big[ II^{1+i\tau}_{l,1} F(x)+  II^{1+i\tau}_{l,2} F(x)\big]
\]
where
\begin{align*}
II^{1+i\tau}_{l,1} F(x)&=
  \intslashx\Big(\sum_{2^k\diam(Q_x )> 2^{10 l}}
\big|
\mathcal A^{k,l} [\chi_{\mathcal U(x)} f_k](y)\big |^2\Big)^{1/2} dy,
\\
II^{1+i\tau}_{l,2} F(x)&=
  \intslashx\Big(\sum_{2^k\diam(Q_x )> 2^{10 l}}
\big|
\mathcal A^{k,l}
 [\chi_{\mathbb R^3\setminus \mathcal U(x)} f_k](y)\big |^2\Big)^{1/2} dy.
\end{align*}
The term $II^{1+i\tau}_{l,1} F(x)$ is estimated by an  $L^2$ estimate; we obtain
after applying the  Cauchy-Schwarz  inequality, (almost) orthogonality of the
$\mathcal A^{k,l}$ and (\ref{L-2}),
\begin{align*}
|II^{1+i\tau}_{l,1} F(x)|&\le
\Big(  \frac{1}{|Q_x |}\int \sum_{2^k\diam(Q_x )> 2^{10 l}}
\big|
\mathcal A^{k,l} [\chi_{\mathcal U(x)} f_k](y)\big |^2 dy \Big)^{1/2},
\\
&\lesssim \sup_{2^k\diam(Q_x )> 2^{10 l}}\|\mathcal A^{k,l}\|_{L^2\to L^2}
\Big(  \frac{1}{|Q_x |}\Big \|\Big(  \sum_{k}
|\chi_{\mathcal U(x)}f_k |^2\Big)^{1/2}\Big\|_2^2\Big)^{1/2}
\\
&\leq C_{\epsilon} \sup_{2^k\diam(Q_x )> 2^{10 l}}2^{(l-k)/2 + l\epsilon}
 \Big( \frac{|\mathcal U(x)|}{|Q_x |}\Big)^{1/2}\|F\|_{L^\infty(\ell^2)}
\\
&\leq C_{\epsilon}  2^{-9 l/2 + l\epsilon}\|F\|_{L^\infty(\ell^2)},
\end{align*}
where at the last step we have used the fact that $|\mathcal
U(x)|/|Q_x| \lesssim \diam(Q_x)^{-1}$.

We now crudely estimate  the terms
$II^z_{l,2} F(x)$ and $III^z_l F(x)$, $z=1+i\tau$.
For this we make use of the following pointwise estimate obtained from
integration by parts :
\begin{multline}
|\mathcal A^{k,l}  f(y)|+2^{-k} |\nabla \cA^{k,l}  f(y)|
\\
\lesssim \iint \frac {2^{3k -2l}}{(1+ 2^{k-2l}|y-w+\gamma(s)|)^N}|f(w)| dw \,ds,
\label{pointwise}
\end{multline}
with $N\geq 4$. To estimate
$|II^{1+i\tau}_{l,2} F(x)|$ we need \eqref{pointwise} for $y\in Q_x $
and $w\notin
\mathcal U(x)$, i.e. $|y-w+\gamma(s)|\gtrsim \diam(Q_x )$ for all relevant $s$.
This yields the bound
\begin{align*}
|T^{1+i\tau}_{l,k}f_k(x)|
&\lesssim
\int_s\frac{1}{ |Q_x|}
\int_{Q_x}
\int_{\substack{
|y-z+\gamma(s)|\\ \ge c \diam(Q_x )  }}
\frac{ 2^{3k-2l} 2^{(2l-k)N}}{|y-z+\gamma(s)|^N}
|f_k(z)| dz \,dy \,ds 
\\
&\lesssim
2^{3k-2l} 2^{(2l-k)N} (\diam (Q_x) )^{3-N} \|F\|_{\ell^\infty(L^\infty)}.
\end{align*}
Therefore,
\[
\sum_{2^k\diam(Q_x )> 2^{10 l}}
|T^{1+i\tau}_{l,k}f_k(x)|
\lesssim 2^{2l(N-1)} 2^{10 l (3-N)}  \|F\|_{\ell^\infty(L^\infty)},
\]
which certainly implies
\begin{equation}
 |II^{1+i\tau}_{l,2} F(x)|\lesssim 2^{-4l} \|F\|_{L^\infty(\ell^2)}. \label{II-infty}
\end{equation}
To estimate $III^z_l F(x)$  we use instead the estimate for the
gradient in (\ref{pointwise}) and get
\begin{align*}
|T^{1+i\tau}_{l,k}f_k(x)| &\le \int_{Q_x}\int_{Q_x}
%\intslashx\intslashx
\Big| \int_{\sigma=0}^1
\langle{y-z},{\nabla \mathcal A^{k,l} f(y+\sigma(z-y))}\rangle d\sigma\Big|
\frac{dz}{|Q_x|}
\frac{dy}{|Q_x|}
%dz dy
\\&\le
2^{4l} 2^k \diam (Q_x ) \|f_k\|_\infty.
\end{align*}
We sum over $k$ with $2^k\diam(Q_x )\le 2^{-10 l}$ and obtain
\begin{equation}
|III^{1+i\tau}_l F(x)|
\lesssim 2^{-6l} \|F\|_{\ell^\infty(L^\infty)}
\lesssim 2^{-6l} \|F\|_{L^\infty(\ell^2)}. \label{III-infty}
\end{equation}
Interpolating the bounds (\ref{II-infty}) and (\ref{III-infty}) with (\ref{II-two}) we obtain
(\ref{maincz}) with $\epsilon(p)>0$ for a range of $p$'s  which includes $(4,\infty)$ and
therefore $((p_W+2)/2, \infty)$.

We observe that by choosing a parameter larger  than $10$ in the definition
of $I,II,III$  we could enlarge the range where $\epsilon(p)>0$ in (\ref{maincz}), but this
is irrelevant here.
This  finishes the proof of Theorem \ref{fixedt} in the case of
nonvanishing curvature and torsion. \qed

\medskip

{\noi \bf Extension to finite type curves.}

We now consider the averaging operator $\cA_t$
%%, $1\le t\le 2$
 as in \eqref{At} and assume that $\gamma$ is of maximal type $n$.
We shall fix $s_0$ and estimate
$\cA_t$ under the assumption that the cutoff function $\chi$
is supported in a small neighborhood of $s_0$.
This  assumption implies that there are orthogonal unit vectors
$\theta_1,\theta_2,\theta_3$ and integers $1\le n_1<n_2<n_3\le n$ so that
for $i=1,2,3$,
$$
\inn{\theta_i}{\gamma^{(j)}(s_0)}=0,
\text{ if $1\le j<n_i$},\quad
\inn{\theta_i}{\gamma^{(n_i)}(s_0)}\neq 0.
$$
After a rotation we may also assume that $\theta_i=e_i$, $i=1,2,3$,
%By a translation we may assume that $\gamma(s_0)=0$ and by
% reparametrization we may assume $s_0=0$,
and
\begin{equation}\label{expansionofgamma}\gamma(s_0+ \a)=\gamma(s_0)+ (\beta_1 \a^{n_1}\varphi_1(\a),
\beta_2 \a^{n_2}\varphi_2(\a),
\beta_3 \a^{n_3}\varphi_3(\a))
\end{equation}
where $\beta_1$, $\beta_2$, $\beta_3$ are nonzero constants and
$\varphi_i\in C^{n+5-n_i}$ with
$\varphi_i(0)=1$.
Thus we need to establish the asserted $L^p\to L^p_{1/p}$-boundedness for
 the averages
$$\cA_t f(x) =\int \chi(\a)
f(x-t\gamma(s_0+\a)) d\a
$$
with bounds uniformly in $t\in [1/2,2]$, where $\chi$ is chosen so that  we assume that $1/2\le |\varphi_i(\a)|\le 3/2$, $i=1,2,3$,  in the support of $\chi$.
We work with a dyadic partition of unity $\zeta_j$, where
$\zeta_j=\zeta(2^j\cdot)$ is supported where $|\a|\approx 2^{-j}$;
we also set $\chi_j=\chi\zeta_j(2^{-j})$ so that the derivatives of
$\chi_j$ are  bounded  independently of $j$.
Let
\begin{align}
A_{j,t} f(x):&=
\int \chi_j(2^j\a)
f(x-t\gamma(s_0+\a)) d\a
\notag
\\
\label{definitionofAtj}
&= 2^{-j}
\int \chi_j(u)
f(x-t\gamma(s_0+2^{-j}u)) du
\end{align}
so that $\cA_t=\sum_{j>0} A_{j,t}$.
Now set
\begin{equation}\label{finitetypedilation}
\begin{aligned}
\delta_j(x)&=(2^{jn_1}x_1,2^{jn_2}x_2,2^{jn_3}x_3)
\\
\Gamma_j(u)&=
\big(\beta_1u^{n_1} \varphi_1(2^{-j} u),
\beta_2 u^{n_2} \varphi_2(2^{-j} u),
\beta_3 u^{n_3} \varphi_3(2^{-j} u)\big)
\end{aligned}
\end{equation}
A change of variable shows that
\begin{multline}\label{finitetypesubstitution}
A_{j,t} f(x)\,=\, \\ 2^{-j} \int f_{-j}(\delta_jx -t\delta_j\gamma(s_0)-t\Gamma_j(u))
\chi_j(u) du, \,\text{ where } f_{-j}(y)=f(\delta_{-j}y).
\end{multline}
We note that the curves $\Gamma_j$ have $C^{n+5-n_j}$ bounds
(in particular $C^5$ bounds)
independent of $j$, and that the parameter $u$ belongs to   the union of two intervals $\pm (c_1,c_2)$ away from the origin
(with $c_1,c_2$ independent of $j$).
Moreover
$$
\det\begin{pmatrix} \Gamma_j'(u)&\Gamma_j''(u)&\Gamma_j'''(u)
\end{pmatrix} \approx \beta_1\beta_2\beta_3 u^{n_1+n_2+n_3}(1+O(2^{-j}))
$$
so that the uniform results
in the case of nonvanishing curvature
and torsion apply.
Observe that for $\alpha\ge 0$
$$\|g\circ\delta_j\|_{L^p_\alpha}\lc 2^{j(n_3\alpha  -N/p)}
\|g\|_{L^p_\alpha},\qquad
 N=n_1+n_2+n_3.$$ Thus  for
$p > (p_W + 2)/2$,
\begin{align*}\big\|A_{j,t}f
\big\|_{L^{p}_{1/p} }&\lc 2^{-j}  2^{j(n_3/p -N/p)}
\Big\|
 \int f_{-j}(\cdot  -t\delta_j\gamma(s_0)-t\Gamma_j(u)) \chi_j(u) du\Big\|_{L^p_{1/p}}
\\&\lc
2^{j(-1+n_3/p -N/p)}\|f\circ\delta_{-j}\|_p=
2^{j(-1+n_3/p)}\|f\|_p.
\end{align*}
Since we  assume that $p>n\ge n_3$
 we can sum in $j$ to arrive at the desired conclusion.\qed

\section{Microlocal  Smoothing Estimates  for Curves in $\bbR^d$}
\label{LSM1}

In this section we consider a
$C^3$ curve
$u\mapsto \Gamma(u)$ in $\bbR^d$, defined in a compact interval $J$, and we assume that there is a constant $B\ge 1$ so that
$B^{-1}\le |J|\le B$ and for all
$u\in J$
\begin{equation}\label{generalassumptionsonGamma}
 \sum_{i=1}^3 |\Gamma^{(i)}(u)|\le B.
\end{equation}

We study the  space-time  smoothing properties  of the averaging operator,
when localized to the region where
$|\inn{\Gamma''(u)}{\xi}|\approx |\xi|.$
Consider for a compactly supported symbol $a$  the operator
defined by
\begin{equation}\label{definitionAaf}
\fA_\Gamma [a, f] (x,t) = (2\pi)^{-d} \iint a(u,t, \xi)
e^{i\inn{x}{\xi} -it \inn{\Gamma(u)}{\xi}}
\widehat f(\xi) \,du\,d\xi.
\end{equation}

\begin{thm}\label{genmicls}
Let $J_0$ be the closed  subinterval of $J$
with same center and length $|J|/2$, and assume that
 $a$ is supported in $$J_0\times [1/2,2]\times \{\xi:2^{k-1}\le |\xi|
\le 2^{k+1}\}$$ and that the inequalities
\begin{equation}\label{aassumption}
|\partial_u^{i_1} \partial_t^{i_2}\partial_\xi^\alpha a(u,t,\xi)|
\le \cC[a] |\xi|^{-|\alpha|},
\end{equation} hold  for $|\alpha|\le d+2$,
$0\le i_1\le 1$, $0\le i_2\le 1$.
Moreover assume that
\begin{equation}\label{Gammasecondlowerbound}
\big|\inn{\Gamma'(u)}{\xi}\big |+\big|\inn{\Gamma''(u)}{\xi}\big |
\ge B^{-1}|\xi| \quad \text{ if  } (u,t,\xi)\in \supp(a).
\end{equation}
Then for $p>p_W$ and $f\in L^p(\bbR^d)$
\begin{equation}\label{Aafinequality}
\big\|\fA_\Gamma[a,f]\big\|_{L^p(\bbR^{d+1})} \le C(\eps,p, B,d)\, \cC[a]\,
2^{-k (\frac 2p-\eps)}  \|f\|_{L^p(\bbR^{d})}.
\end{equation}
\end{thm}

The crucial hypothesis on $\Gamma$ is the lower bound
\eqref{Gammasecondlowerbound}.
We note that the derivatives  of $\Gamma$ are assumed to be bounded but
we make no size assumption on $|\Gamma(u)|$ itself.
Thus the assumptions on $\Gamma$
are invariant under translation of the curve.

In the following subsection we shall prove this theorem under
slightly more restrictive normalization assumptions which will
be removed at the end of this section by localization and scaling
arguments.

\medskip

\subsection{Normalization}
We now work with a $C^3$ curve $s\mapsto \gamma(s)$, $\gamma:I^*\to
\bbR^d$, $d\ge 3$, where $I^*=[-2\delta,2\delta]$ is a closed subinterval
of $[-1,1]$;
we also set  $I=[-\delta,\delta]$.
We assume that the curve is parametrized by
arclength, i.e. $|\gamma'(s)|=1$ for all $s\in I^*$ and that for
some $M\ge 10,$
\begin{equation}\label{gammaassumptionsRd}
 |\gamma''(s)|+|\gamma'''(s)|\le M.
\end{equation}
Let $\Omega$ be  an open convex conic subset of $\bbR^d\setminus \{0\}$, and let
\begin{equation}\label{definitionofOmegak}
\Omega_k=\{\xi\in \Omega: 2^{k-1}\le |\xi|\le 2^{k+1}\}.
\end{equation}
We shall study $\fA_\gamma[b,f]$ defined as in
\eqref{definitionAaf} and we
%but will omit the subscript $\gamma$ in what follows. We
now assume that the symbol $b$ is supported in $I\times[1/2,2]\times \Omega_k$ and
satisfies
\begin{equation}\label{bassumption}
|\partial_s^{i_1} \partial_t^{i_2}\partial_\xi^\alpha b(s,t,\xi)|
\le C_{i_1,i_2,\alpha} |\xi|^{-|\alpha|},
\end{equation} for $|\alpha|\le d+2$,
$0\le i_1\le 1$, $0\le i_2\le 1$.
We assume that  $\Omega$ satisfies the   crucial

\medskip

\noi{\it Nondegeneracy Hypothesis:}
\begin{equation}\label{gamma2lowerbound}
  |\xi|/2\le |\inn{\gamma''(s)}{\xi}| \le 2|\xi|
\end{equation}
for all $\xi\in \Omega$, $s\in I$; moreover we assume that
for every  $\xi$ in $\Omega$  there is at least one
$s\in [-3\delta /4,3\delta /4]$ so that
\begin{equation}\label{gamma1xismall} |\inn{\gamma'(s)}{\xi}|\le \delta |\xi|/10.
\end{equation}

\medskip

Note that the smallness assumption
\eqref{gamma1xismall} and the lower bound
\eqref{gamma2lowerbound} imply that  for each $\xi\in \Omega$
there is a unique $s=\scr(\xi)$ in $(-\delta ,\delta )$
so that
\begin{equation}\label{implicitfunction}
\inn{\gamma'(s)}{\xi}=0 \,\iff\,  s=\scr(\xi),
\end{equation}
and   $\xi\mapsto \scr(\xi)$
is a $C^3$ function on $\Omega$
which is homogeneous of degree $0$.

The next  subsection is devoted to the proof of
\begin{thm}\label{micls}
Assume that $k\ge 10$, $\eps_0>0$, that
 $b$ is supported in $I\times [1/2,2]\times \Omega_k$ and
that
\eqref{gammaassumptionsRd}
\eqref{bassumption},
and the nondegeneracy hypothesis hold.

Then for $p>p_W$
\begin{equation}\label{Abfinequality}
\big\|\fA_\gamma[b,f]\big\|_{L^p(\bbR^{d+1})}
\le C(\eps_0,p,M,d) 2^{-k (\frac 2p-\eps_0)}
 \|f\|_{L^p(\bbR^{d})},\quad \eps_0>0.
\end{equation}
\end{thm}

\medskip

\subsection{Proof of Theorem \ref{micls}} We use the following

\medskip

\noi{\it Notation:} ``Constants'' $C$ may depend on $M$ and the dimension;
we shall use the Landau symbol $\cE=O(B)$ if $|\cE|\le CB$.
We shall also use the
notation
  $\cE=\cO_1(B)$ if  $|\cE|\le B$.
%(this is the Landau symbol with constant $1$).

\medskip

\noi {\it Some symbol classes}.
 Let $2^{-k/2}<r\le \min\{10^{-3} M^{-2},\delta/4\}$, $s\in (-3\delta/4,3\delta/4)$.
We define some symbol classes for multipliers
$m(\xi,\tau)$ and set $\Xi=(\xi,\tau)$, with
$\tau\equiv \Xi_{d+1}$. We denote by $e_1$, ...,
$e_{d+1}$ the standard basis in
 $\bbR^{d+1}$.
Let $L^{(1)}_{s}$ be the linear shear transformation which maps
$$\Xi=\sum_{i=1}^d\xi_ie_i +\tau e_{d+1}
%\quad \text{  to } \quad
\,\,\mapsto \,\,
L^{(1)}_{s}\Xi=
\sum_{i=1}^d\xi_ie_i +(\tau-\inn{\gamma(s)}{\xi}) e_{d+1}.$$
Let $L^{(2)}_{r,s}$ be the dilation which satisfies
\begin{equation*}
\begin{gathered}
L^{(2)}_{r,s} e_{d+1}= r^{2} e_{d+1}, \quad L^{(2)}_{r,s}
\gamma'(s)=r\gamma'(s)
\\
L^{(2)}_{r,s}v=v \text{ if } v\in \big(\text{span}\{e_{d+1},\gamma'(s)\}
\big)^{\perp};
\end{gathered}
\end{equation*}
here we identify with a slight abuse of notation the function  $\gamma$
with the function $s\to (\gamma(s),0)$ with values in $\bbR^{d+1}$.
We  define the composition
$$L_{r,s}= L^{(1)}_sL^{(2)}_{r,s}.$$
%%and $L^2_{r,s}v=v$ for any vector which is orthogonal to both
%%$e_{d+1}$ and $\gamma'(s)$. The linear transformation $L_{r,s}$ is
%%defined as the composition $L^1_s\circ L^2_{r,s}$.
%%

Let  $\cS_k(r,s)$ be the class of multipliers $m(\xi,\tau)$
which are supported in
\begin{multline}\label{supportassumptionsonSkrs}
\big\{\Xi=(\xi,\tau):
2^{k-1}\le |\xi|\le 2^{k+1}, \quad
\\
|\inn{\gamma'(s)}{\xi}|\le 2^{k+3}r,
\quad|\tau+\inn{\gamma(s)}{\xi}|\le 2^{k+4} r^2\big\}
\end{multline}
%| L_{r,s}^{-1}\Xi| \le 2^{k+3}\big\}
and satisfy
\begin{equation}\label{Hormandermult}
\big| \partial_{\Xi}^\alpha\big(
m(L_{r,s}\Xi)\big)\big|\le |\Xi|^{-\alpha}, \quad |\alpha|\le d+2.
\end{equation}
Note that if  $m\in \cS_k(r,s)$ and $(\xi,\tau)\in\supp (m)$ then
%$|\tau+\inn{\gamma(s)}{\xi}|\le 2^{k+1}r^2$
%and $|\inn{\gamma'(s)}{\xi}|\le 2^{k+1}r$.
$| L_{r,s}^{-1}\Xi| =O(2^{k})$.
The following Lemma is straightforward to check, we omit the proof.
\begin{lemmasub}\label{varyings}
 There
are  constants $C_i=C_i(A,M)$, $i=1,2$,
so that $C_1^{-1} m\in \cS_k(C_2r,s')$
 for all
$m\in \cS_k(r,s)$ and
all $s'$ with
$|s-s'|\le A r$.
\end{lemmasub}

We shall need kernel estimates for operators associated with
 multipliers in $\cS_k(r,s)$.

\begin{lemmasub}
\label{lemmapointwiseKsestimate}
 Let $m(t,\cdot) \in \cS_k(r,s')$ for $1/2\le t\le 2$
(depending continuously on $t$)
and assume $|s'-s|\le 2r$.
Let
\begin{equation}\label{Ksm}
K_s[m](x,t')=\iiint e^{i\inn{x}{\xi}+it'\tau}
e^{-it(\tau+\inn{\gamma(s)}{\xi})}m(t,\xi,\tau)
\, d\xi d\tau \,dt.
\end{equation}
Then
\begin{multline}\label{pointwiseKsestimate}
|K_s[m](x,t')|\le C\, \times
% C \int
%(1+2^k r^2|t-t'|)^{-2}\quad C \times
\\
 \int
\frac{2^{k(d+1)}r^3}
{(1+2^k r^2|t-t'|+  2^k r|\inn{x-t'\gamma(s)}{\gamma'(s)}|
+2^k |\Pi^\perp_{\gamma'(s)}(x-t'\gamma(s))|)^{d+2}}\, dt
\end{multline}
where $\Pi^\perp_{\gamma'(s)}:\bbR^d\to \bbR^d$ denotes the orthogonal projection to the orthogonal complement of $\bbR \gamma'(s)$.
In particular
$$
\sup_{x,t'}\int|K_s[m](x-y,t')| dy \le C.
$$
\end{lemmasub}

\begin{proof}
The second assertion is an immediate consequence of
 \eqref{pointwiseKsestimate}.
To see
\eqref{pointwiseKsestimate} we
 change variables in the
integral defining  $K_s$ and see that
\begin{equation}\label{Ksmafterchange}
K_s[m](x,t')=
\int e^{i(t'-t)\tau}
e^{i\inn{x-t'\gamma(s)}{\xi}}
m(t,L^{(1)}_s\Xi)  d\Xi dt
\end{equation}
with $\Xi=(\xi,\tau)$.
Changing variables again using $\Xi=L^{(2)}_{r,s}\widetilde \Xi$ and
integrating
by parts in  $\widetilde \Xi$ yields \eqref{pointwiseKsestimate}.
\end{proof}

Now let $s_0\in (-3\delta/4, 3\delta/4)$, $r\le \delta/8$ and let
$\fS_k(r,s_0)$ be the class of symbols
$(s,t,\xi,\tau)\mapsto a(s,t,\xi,\tau)$ which are supported in
$[s_0-2r, s_0+2r]\times[1/2,2]\times \Omega_k\times \bbR$
%$2^{k-1}\le |\xi|\le 2^{k+1}$,  $|L_{r,s}^{-1} (\Xi)|\le 2^{k+1}$,
%$|\inn{\gamma''(s)}{\xi}|\ge C_1 2^k$, $|s-s^0|\le 2r$, $|t|\le 2$
and which satisfy
\begin{equation*}
\partial_s^{i_1}\partial_t^{i_2}
 a(s,t,\cdot) \in \cS_k(r, s_0),
\quad\text{
for
$i_1, i_2 \in \{0,1\}$.}
\end{equation*}

%Given such a symbol
We  define the oscillatory integral
\begin{equation}\label{defTaf}
\cT[a, f](x,t)=(2\pi)^{-d-1}\iint e^{i(\inn{x}{\xi}+t\tau)}
\fm[a](\xi,\tau)
\widehat f(\xi) d\xi d\tau
\end{equation}
where
\begin{equation}\label{defma}
\fm[a](\xi,\tau)=\iint e^{-it(\tau+\inn{\gamma(s)}{\xi})} a(s,t,\xi,\tau) ds dt.
\end{equation}

The proof of Theorem  \ref{micls} relies on an iteration where
the main step is to prove  the following proposition. Here we say that a set
of real numbers  is {\it $r$-separated} if $|s-s'|\ge r$ for different
$s$, $s'$ in this set.

\begin{propositionsub}
\label{maininductionstep}
 Suppose that  $\eps>0$,
$r_1\le r_0\le \min\{10^{-3}M^{-2}, \delta/4\}$,
and assume that
$r_1\ge 100 Mr_0^{3/2}$.
 Let
 $\{s^\mu\}$ be an $r_0$-separated set of points  in $[-\delta,\delta]$,
 and for each $\mu$
let
$a^\mu$ be a  symbol in $\fS_k(r_0, s^\mu)$. Let $p>p_W$. Then
there is
a set of $r_1$-separated points $\{s_\nu\}$ and symbols $a_\nu\in
\fS_k(r_1, s_\nu)$, and for every $\eps>0$ there is a constant
$\cC_\eps=\cC_\eps(p,M)$,
so that
\begin{multline}\label{propositionin4}
\Big(\sum_\mu\big\| \cT[a^\mu,f]\big\|_p^p\Big)^{1/p}
\\ \le \cC_\eps (r_0/r_1)^{1-\frac 4p +\eps}
\Big[
\Big(\sum_\nu\big\| \cT[a_\nu,f]\big\|_p^p\Big)^{1/p} +
2^{-k}r_1^{-1}\|f\|_p\Big]
\end{multline}
holds.
\end{propositionsub}

\noi {\bf Proof of the Proposition.}

For each $\mu$ we set
\begin{equation}\label{Umu}
U_\mu(\xi,\tau) = \tau+\inn{\gamma(s^\mu)}{\xi}-
\frac {1}{2} \frac{\inn{\gamma'(s^\mu)}{\xi}^2}
{\inn {\gamma''(s^\mu)}{\xi}}.
\end{equation}

The second part of the  following lemma states
%%that a localization in $s$ corresponds to a
%%localization in terms of
%%$\inn{\gamma'(s^\mu)}{\xi}/
%%\inn{\gamma''(s^\mu)}{\xi}$ and
that $U_\mu$ is a good approximation for
$\tau+\inn{\gamma(\scr(\xi))}{\xi}$.

\begin{lemmasub}\label{estimatesonsuppamu}
Suppose $|s|\le 3\delta/4$, $\xi\in \Omega$.
%Suppose $(s,\xi,\tau)$
%belongs to the support of
%$a^\mu$.
 Then
\begin{equation}
\label{approximationone}
s-\scr(\xi)
=\frac{\inn{\gamma'(s)}{\xi}}
{\inn{\gamma''(s)}{\xi}} +\cO_1(6M(s-\scr(\xi))^2);
\end{equation}
in particular this holds for $s=s^\mu$ if  $(\xi,\tau) \in \supp (a^\mu)$
for some $\tau$.
Moreover
\begin{equation}
\label{approximationtwo}
U_\mu(\xi,\tau)\,=\,\tau+\inn{\gamma(\scr(\xi))}{\xi}+
\cO_1(13M |\scr(\xi)-s^\mu|^3)|\xi|
\end{equation}
if $(\xi,\tau) \in \supp (a^\mu)$.
\end{lemmasub}
\begin{proof}
We expand using \eqref{implicitfunction}
\begin{align}
\inn{\gamma'(s)}{\xi}&= \inn{\gamma''(\scr(\xi))}{\xi} (s-\scr(\xi))+\cO_1
(M|\xi| (s-\scr(\xi))^2/2)
\notag
\\
&= \inn{\gamma''(s)}{\xi}(s-\scr(\xi))
+\cO_1
(3M|\xi| (s-\scr(\xi))^2/2)
\label{approximationthree}
\end{align}
and \eqref{approximationone}
follows by using the lower bound in
\eqref{gamma2lowerbound}.

Next expand again using \eqref{implicitfunction}
\begin{align*}
%% \label{approximationfour}
\inn{\gamma(s^\mu)}{\xi}&-\inn{\gamma(\scr(\xi))}{\xi}
\\=&
\inn{\gamma''(\scr(\xi))}{\xi}\frac{(s^\mu-\scr(\xi))^2}2 +
\cO_1(M|\xi| |s^\mu-\scr(\xi)|^3/6)
\\
=&
\inn{\gamma''(s^\mu)}{\xi}\frac{(s^\mu-\scr(\xi))^2}2 +
\cO_1(2M|\xi| |s^\mu-\scr(\xi)|^3/3).
\end{align*}
Now we use  \eqref{approximationone}
for $s=s^\mu$ and get
\begin{align*}
&\inn{\gamma''(s^\mu)}{\xi}\frac{(s^\mu-\scr(\xi))^2}2 -
\frac{\inn{\gamma'(s^\mu)}{\xi}^2}
{2\inn{\gamma''(s^\mu)}{\xi}}
\\&= \inn{\gamma'(s^\mu)}{\xi} \cO_1(6M(s^\mu-\scr(\xi))^2)
+\cO_1(18 M^2(s^\mu-\scr)^4)|\xi|
\\
&=
 \cO_1(12 M(s^\mu-\scr(\xi))^3)|\xi|+\cO_1(20 M^2(s^\mu-\scr)^4)|\xi|.
\end{align*}
%In the last step we have used the upper bound on $\gamma''$ in
%\eqref{gammaassumptionsRd}.
 Since we assume that $r\le 10^{-3}M^{-1}$
we obtain \eqref{approximationtwo}.
\end{proof}

We now decompose $a^\mu$ using cutoff functions $\eta_0$, $\eta_1$, $\zeta$
as in  \S\ref{averages}, that is,
$\eta_0$ is supported in $[-1,1]$, equal to $1$ in $[-1/2,1/2]$,
$\eta_1=\eta_0(4^{-1}\cdot)-\eta_0$, and $\zeta$ is supported
in $(-1,1)$ and satisfies $\sum_\nu \zeta(s-\nu)=1$, $s\in\bbR$. These cutoff functions are fixed and various constants below
may depend on their choice.
Set
\begin{equation*}
a^\mu_{0,\nu}(s,t,\xi,\tau)
=a^\mu(s,t,\xi,\tau)
\eta_0\big(r_1^{-2}(2^{-k}|U_\mu(\xi,\tau)|+(s-\scr(\xi))^2)\big)
 \zeta(r_1^{-1}s-\nu),
%\big(A \tfrac{(s-\scr(\xi))^2}{U_\mu(\xi,\tau)}\big)
%\eta_0 (r_1^{-2} U_\mu(\xi,\tau) \zeta(r_1^{-1}s-\nu),
\end{equation*}
and, for $n\ge 1$

%\begin{equation*}
\begin{align*}
a^\mu_{n,\nu}(s,t,\xi,\tau)
=&a^\mu(s,t,\xi,\tau)
\eta_1\big(2^{2-2n} r_1^{-2} (2^{-k}|U_\mu(\xi,\tau)|+(s-\scr(\xi))^2)\big)
\\& \, \times\,
\eta_0 \big(\frac{(s-\scr(\xi))^2}{2^{-k-8}U_\mu(\xi,\tau)}\big)
\zeta(2^{-n}r_1^{-1}s-\nu),
%\quad n\ge 1,
%\eta_1 (2^{2-2n}r_1^{-2} U_\mu(\xi,\tau))
%\zeta(2^{-n}r_1^{-1}s-\nu), \quad n\ge 1,
\\
b^\mu_{n,\nu}(s,t,\xi,\tau)
=&a^\mu(s,t,\xi,\tau) \, \eta_1\big(2^{2-2n}r_1^{-2}(2^{-k}|U_\mu(\xi,\tau)|+(s-\scr(\xi))^2)\big)
\\& \times\,
\big(1-\eta_0 \big(\frac{(s-\scr(\xi))^2}{2^{-k-8}U_\mu(\xi,\tau)}\big)\big)
 \zeta(2^{-n}r_1^{-1}s-\nu)
%\eta_0(r_1^{-1}(s-\scr(\xi)) \zeta(r_1^{-1}s-\nu).
\end{align*}
Then
\begin{equation}\label{amudecomposition}
a^\mu =\sum_\nu a^{\mu}_{0,\nu}+
\sum_{n\ge 1} \sum_\nu(a^{\mu}_{n,\nu} + b^{\mu}_{n,\nu}).
\end{equation}
%We observe that
%$a^\mu_{0,\nu}$ is supported where
%$|U_\mu(\xi,\tau)|\lc 2^k r_1^2$,
%$|s-\scr(\xi)|\lc r_1$ and $|r_1\nu-s|\lc r_1$.
%Moreover, for $n\ge 1$,
%$a^\mu_{n,\nu}$ is supported where
%$|U_\mu(\xi,\tau)|\approx 2^k 2^{2n} r_1^2$, $|s-\scr(\xi)|\lc
%2^n r_1$, $|s-\scr(\xi)|^2\ll 2^{-k} U_\mu$
% and $|r_1\nu-s| \lc 2^n r_1$. Finally,
%$b^\mu_{n,\nu}$ is supported where
%$|s-\scr(\xi)|\approx  2^{n} r_1$,
%$|U_\mu(\xi,\tau)| \lc 2^{k+2n} r_1^2$,
%$2^{-k}|U_\mu(\xi,\tau)| \lc (s-\scr(\xi))^2$,
%$|r_1\nu-s|\lc 2^n r_1$.

In what follows we define
the linear map $\omega^\mu
%=(\omega^\mu_1, \omega^\mu_2, \omega^\mu_3)
:\bbR^{d+1}\to \bbR^3$ by
%For $(\xi,\tau)\in \bbR^{d+1}$ let
\begin{equation}\label{omegamu}
%\omega^\mu(\xi,\tau)=
%\big(\inn{\gamma'(s^\mu)}{\xi},
%\tau +\inn{\gamma(s^\mu)}{\xi},
%\inn{\gamma''(s^\mu)}{\xi}\big).
\begin{cases}\omega^\mu_1(\xi,\tau)
&=\inn{\gamma'(s^\mu)}{\xi},
\\
\omega^\mu_2(\xi,\tau)
&=
\tau +\inn{\gamma(s^\mu)}{\xi},
\\
\omega^\mu_3(\xi,\tau)
&=\inn{\gamma''(s^\mu)}{\xi}.
\end{cases}
%\omega^\mu_1=\inn{\gamma'(s^\mu)}{\xi},
%\quad
%\omega^\mu_2=
%\tau +\inn{\gamma(s^\mu)}{\xi},\quad
%\omega^\mu_3=\inn{\gamma''(s^\mu)}{\xi}.
\end{equation}
We shall  observe that for fixed $\mu$  the supports   of
$a^{\mu}_{n,\nu}$
and
$b^{\mu}_{n,\nu}$ are contained in ``plates'' defined using
$\omega^\mu(\xi,\tau)$ (\cf.
\eqref{platemunu}
below).

\begin{lemmasub}\label{supportproperties}
%4.2.5
%(i)
%%COMMENT: Now true by definition
%If $(\xi,\tau)$ is in the support of $a^\mu(s,\cdot)$ then
%$|s-s^\mu|\le 2r_0$, and
%%\begin{equation}\label(support1}
%$|\inn{\gamma(s^\mu)}{\xi}|\le 2^{k+1}r_0$ and
%$|\tau+\inn{\gamma(s^\mu)}{\xi}|\le 2^{k+1}r_0^2$.
%%\end{equation}
%
(i) Suppose $(\xi,\tau)$ is in the support of $a^\mu_{n\nu}(s,t,\cdot)$
or $b^\mu_{n\nu}(s,t,\cdot)$
then $|U_\mu(\xi,\tau)| \le 2^{k+2n}r_1^2 $,
and $|s-\scr(\xi)|\le 2^nr_1$.

(ii) Suppose $n\ge 1$.

If $(\xi,\tau)$ is in
 the support of $a^\mu_{n\nu}(s,t,\cdot)$ then
$|U_\mu(\xi,\tau)|\ge 2^{k+2n-3}r_1^2$.

If $(\xi,\tau)$ is in the support of $b^\mu_{n\nu}(s,t,\cdot)$ then
$|s-\scr(\xi)|\ge 2^{n-5}r_1$.

(iii)
Let $s_{n\nu}=2^n r_1\nu$ for $\nu\in \bbZ$ and
assume  that $|s_{n\nu}-s^\mu|\le 2r_0$. Then
there is a constant $C$  so that
the symbols $C^{-1} a^\mu_{n,\nu}$, $C^{-1} b^\mu_{n,\nu}$
belong to $\fS_k(2^n r_1, s_{n\nu})$.

(iv)
If $2^{n}r_1> 2^4 r_0$ then $a^\nu_{n,\nu}=0$.

(v) If $2^n r_1> 2^{7} r_0$ then $b^\nu_{n,\nu}=0$.

(vi)
Let $g(\a) =(\a,\a^2/2)$ and let $u_i(\a)$, $i=1,2,3$, be as in \eqref{uialpha}, i.e.
$$
u_1(\a)= (\a,\a^2/2,1),\  u_2(\a)=(1,\a,0), \ u_3(\a)=(-\a,1,\a^2/2).
$$
Let
$\oa=
\oa_{\mu n\nu}= s^\mu-s_{n\nu}$.
Then the supports of
$a^\mu_{n, \nu}$ and  $b^\mu_{n, \nu}$
are  contained in the set $Pl_{\mu\nu}^n$  consisting of all $(\xi,\tau)$
which satisfy
\begin{equation}
\label{platemunu}
%\{(\xi,\tau):
\begin{cases}
|\inn{u_1(\oa)}{\om^\mu(\xi,\tau)}|\le  2^{k+2},\\
|\inn{u_2(\oa)}{\om^\mu(\xi,\tau)-\omega_3^\mu(\xi,\tau) u_1(\oa)}
|\le  2^{k+4} 2^nr_1,\,
\\
|\inn{u_3(\oa)}{\om^\mu(\xi,\tau)
%%%-\omega_3^\mu(\xi,\tau) u_1(\oa)
}|\le
 2^{k+3} 2^{2n} r_1^2.
%\}
\end{cases}
\end{equation}

(vii) Every $(\xi,\tau)$  belongs to no more than $75$  of the sets
$\{(\xi,\tau): (s,t,\xi,\tau)\in \supp  (a^{\mu}_{n,\nu})\}$ and to no more than $75$ of the sets
$\{(\xi,\tau): (s,t,\xi,\tau)\in \supp  (b^{\mu}_{n,\nu})\}.$
\end{lemmasub}

\begin{proof}
Properties (i) and  (ii) are immediate consequences of the definition
of the symbols.

For (iii) we first have to check the support properties, namely
assuming that $(s,t,\xi,\tau)$ belongs to the support of $a^\mu_{n,\nu}$
or
$b^\mu_{n,\nu}$ then
\begin{align}
\label{firstsuppprop}
|\inn{\gamma'(s_{n\nu})}{\xi}|
&\le 2^{k+3} 2^{n}r_1
\\
\label{secondsuppprop}
|\tau+\inn{\gamma(s_{n\nu})}{\xi}|
&\le
2^{k+4} 2^{2n}r_1^2
\end{align}
To see this we first note that $|s-s_{n\nu}|\le 2^n r_1$
and
$(s-\scr(\xi))^2\le 2^{2n} r_1^2$ hence
\begin{equation}
\label{scrsnnu}
|\scr(\xi)-s_{n\nu}|\le 2^{n+1} r_1\end{equation}
Similarly we have of course also
\begin{equation}
\label{scrsmu}
|\scr(\xi)-s^\mu|\le 2r_0\end{equation}
Now $\inn{\gamma'(s_{n\nu})}{\xi}=
\inn{\gamma''(\tilde s)}{\xi}
(s_{n\nu}-\scr(\xi))$
where $\tilde s$ is between $s_{n\nu}$ and $\scr(\xi)$. Since
$|\inn{\gamma''(\tilde s)}{\xi}|\le 2^{k+2}$ we conclude
\eqref{firstsuppprop}.

To see
 \eqref{secondsuppprop} we expand
$$\tau+\inn{\gamma(s_{n\nu})}{\xi}=
\tau+\inn{\gamma(\scr(\xi))}{\xi}
+\inn{\gamma'' (s')}{\xi}\frac{(s_{n\nu}-\scr(\xi))}{2}
$$ where
$s'$ is between $s_{n\nu}$ and $\scr(\xi)$. From \eqref{approximationtwo}
and  \eqref{scrsnnu}
we obtain
$$
|\tau+\inn{\gamma(s_{n\nu})}{\xi}|\le 2^{k+3} 2^{2n}r_1^2+
|U_\mu(\xi,\tau)|+13 M (2r_0)^3|\xi|.$$
Now $|U_\mu(\xi,\tau)|\le 2^{k+2n}r_1^2$ and from our crucial assumption
on the relation between $r_0$ and $r_1$, namely $r_0^3\le (100 M)^{-2} r_1^2$, we can deduce \eqref{secondsuppprop}.

We now have to verify the symbol estimates \eqref{Hormandermult}.
First observe $\partial_\tau U_\mu=1$ and calculate  (using the
notation in \eqref{omegamu})
\begin{align*}
\nabla_\xi U_\mu=\gamma(s^\mu)-\frac{\omega_1^\mu}{\omega_3^\mu}
\gamma'(s^{\mu})+\frac12 \Big(\frac{\omega_1^\mu}{\omega_3^\mu}
\Big)^2 \gamma''(s^\mu)
\end{align*}
and an expansion about the point $s_{n\nu}$ yields that
%$\nabla_\xi U_\mu=$
%\begin{align*}
\begin{multline*}
\nabla_\xi U_\mu=\gamma(s_{n\nu})+ \gamma'(s_{n\nu})\big(s^\mu-s_{n\nu}-
\frac{\omega^\mu_1}{\omega^\mu_3}\big)\\
+\gamma''(s_{n\nu}) \Big(
\frac{(s^\mu-s_{n\nu})^2}{2} -
\frac{\omega_1^\mu}{\omega_3^\mu}(s^\mu-s_{n\nu})+(\frac 12\Big(
\frac{\omega^\mu_1}{\omega^\mu_3}\Big)^2\Big)+O(r_0^3).
\end{multline*}
A further expansion using
\eqref{approximationthree}
 shows that on the support of either $a^{\mu}_{n,\nu}$ ($n\ge 0$) or
$b^{\mu}_{n,\nu}$ ($n\ge 1$)
\begin{equation*}
%\label{gradientUmuestimate}
\nabla_\xi U_\mu
=\gamma(s_{n\nu})+ \gamma'(s_{n\nu}) O( 2^n r_1)+
O(2^{2n} r_1^2)+O(r_0^3).
\end{equation*}
Thus if $v$ is perpendicular to $(\gamma(s_{n\nu}),1)$ then
$\inn{v}{\nabla} U_\mu=O(2^{n} r_1)$ and if $v$ is perpendicular to
both
$(\gamma(s_{n\nu}),1)$ and $(\gamma'(s_{n\nu}),0)$ then
$\inn{v}{\nabla} U_\mu=O(2^{2n} r_1^2)$.
Moreover, from \eqref{implicitfunction}
$$\nabla_\xi \scr(\xi)=\frac{-\gamma'(\scr(\xi))}
{\inn{\gamma''(\scr(\xi))}{\xi}}
$$
and thus
$\inn{v}{\nabla} \scr(\xi) =O(2^{n-k} r_1)$ if $v$ is perpendicular to
$(\gamma'(s_{n\nu}),0)$. Given the bounds on the directional derivatives
of $\scr$ and $U_\mu$ the verification of
%the symbol estimates
\eqref{Hormandermult} is
straightforward.

Next to see (iv) observe that
$|\omega^\mu_2|\le 2^{k+3} r_0^2$
on the support of $a^\mu$
(and hence on the support of $a^{\mu}_{n\nu}$).
But we also have $|U_\mu|\ge 2^{k+2n-3}  r_1^2$ and
$$(\omega_2^\mu)^2/|\omega^\mu_3|\le 2^{k+3} |s-\scr(\xi)|^2\le
2^{k+3}(2^{-k-8}|U_\mu|)
$$
and thus
$|\omega_2^\mu|\ge |U_\mu|/2 \ge 2^{k+2n-4} r_1^2$. This forces
$2^{n}r_1\le 2^4 r_0$
if the support of $a^\mu_{n\nu}$ is nonempty.

Next consider the  support of $b^{\mu}_{n,\nu}$ where  $|\scr(\xi)-s^\mu|\le
2r_0$ and also $(\scr(\xi)-s^\mu)^2\ge 2^{-k-9}|U_\mu(\xi,\tau)|$;
moreover
$\max \{(s-\scr(\xi))^2,2^{-k}|U_\mu(\xi,\tau)\}\ge 2^{2n-3}r_1^2$.
These conditions imply that $2^{2n-3}r_1^2 \le 2^9 (2r_0)^2$
if the support of $b^\mu_{n,\nu}$ is nonempty  and thus (v) follows.

To see (vi) we set
$$\ob=\ob(\xi)=\scr(\xi)-s_{n\nu}$$
and observe that $|\ob|\le 2^{n+1} r_1$ in the supports of $a^\mu_{n,\nu}$
and
$b^\mu_{n,\nu}$, moreover $|\oa| \le 2r_0$.
By \eqref{approximationone} for $s=s^\mu$,
%\begin{equation*}
%\inn{\ga'(s^\mu)}{\xi}
%= (s^\mu-\scr(\xi)) \inn{\ga''(s^\mu)}{\xi}+
% \cO_1(M(s^\mu-\scr(\xi))^2|\xi|/2)
%\end{equation*}and thus
%\begin{equation}s^\mu-\scr(\xi)=
%\frac {\inn{\ga'(s^\mu)}{\xi}}{\inn{\ga''(s^\mu)}{\xi}}
%+\cO_1(4Mr_0^2)\end{equation}
we can write
\begin{equation}\label{oaexpansion}
\oa=\omega_1^\mu/\omega_3^\mu+\ob+E \quad\text{ with }
E=\cO_1(24Mr_0^2).
\end{equation}
Also  $
\omega_2^\mu(\xi,\tau)=U_\mu(\xi,\tau) +
(\omega^\mu_1)^2/(2\omega_3^\mu)$ and
%and also  $2^{-k}|U_\mu(\xi,\tau)|\le 2^{2n}r_1^2$ in the supports of
%$a^\mu_{n,\nu}$ and
%$b^\mu_{n,\nu}$.
%This implies that
%$$$|\omega_2^\mu|\le  4r_0^2+2^{2n } r_1^2$$ and
\begin{align*}
&\inn{u_2(\oa)}{\om^\mu-\omega_3^\mu u_1(\oa)}
=\omega_1^\mu -\omega_3^\mu\oa+\omega_2^\mu\oa-\omega_3^\mu \oa^3/2
\\
&= \ob(\om^\mu_2-\omega^\mu_3)+ \big[ E(\omega_3^\mu-\omega_2^\mu)+
\frac{\omega_2^\mu\omega_1^\mu}{\omega_3^\mu}-\omega_3^\mu \oa^3/2
\big]
\end{align*}
and the expression $[...]$ is easily seen to be $\cO_1(100 M r_0^2)$
in view of the assumptions on the support of $a^\mu$. Since we also
assume
$r_1>100 M r_0^{3/2}$ we deduce
$$|\inn{u_2(\oa)}{\om^\mu-\omega_3^\mu u_1(\oa)}|\le 2^{k+n+4} r_1.$$
Next we compute using \eqref{oaexpansion}
\begin{align*}
%\label{u3oa}
\inn{u_3(\oa)}{\om^\mu}&=-\om_1^\mu \oa+
\om_3^\mu \frac{\oa^2}{2}
+\frac{(\om_1^\mu)^2}{2\omega_3^\mu}+U_\mu
\\
&= U_\mu+ \ob^2\om_3^\mu/2+ \omega_1^\mu E+ \omega_3^\mu\ob E+
\omega_3^\mu E^2/2
\\&= \cO_1(2^{k+1} 2^{2n} r_1^2)+\cO_1(100 M 2^kr_0^3).
\end{align*}
which concludes the proof of (vi).

Now suppose $(\xi, \tau) $ belongs to the support of
$a^\mu_{n,\nu}$.
We have seen  that then $|\scr(\xi)-s^\mu|\le 2r_0.$
Since we assume that $\{s^\mu\}$ is a $r_0$-separated set we see that
$(\xi,\tau)\in \supp (a^\mu_{n,\nu})$ for some $n,\nu$ can only happen for
at most five  $\mu$.
Exacly the same argument works for $b^{\mu}_{n\nu}$ in place of
$a^{\mu}_{n\nu}$.
Next, note that for a  $\sigma \in\bbR$ one has
$\eta_1(2^{2-2n}\sigma)\neq 0$ for at most five values of $n$. The
 definition of the   functions $a^\mu_{n,\nu}$ and
$b^\mu_{n,\nu}$ shows that  for fixed $\mu$ there are at most
five values of $n$ for which $a^\mu_{n,\nu}\neq 0$  or
$b^\mu_{n,\nu}\neq 0$.
Finally, given $u$
 there are at most three values of $\nu$ for which $\zeta (u-\nu)\neq 0$.
This shows that for fixed $\mu$ and fixed $n$ there are at most three values of $\nu$ for which
 $a^\mu_{n,\nu}\neq 0$  or
$b^\mu_{n,\nu}\neq 0$. A combination of these observations
  yields the assertion (vii).
\end{proof}

\begin{lemmasub}\label{wolffestimateinRd}
Suppose $Pl^n_{\mu\nu}$
is as in \eqref{platemunu}, and $\mu$, $n$ are fixed.
Suppose that the Fourier transform of $f_\nu$ is supported in
$Pl^n_{\mu\nu}$. Let $\fJ^\mu_n=\{\nu:|s_{n\nu}-s^\mu|\le 2r_0\}.$
Then for $\eps>0$, $p>p_W$,
$$\Big\|\sum_{\nu\in \fJ^\mu_n} f_\nu\Big\|_p\le \cC_{p,\eps}
\Big(\frac{r_0}{2^n r_1}\Big)^{1-\frac 4p+\eps}
\Big(\sum_\nu\|f_\nu\|_p^p\Big)^{1/p}.
$$
\end{lemmasub}

\begin{proof}
Note that $\om_\mu$ in \eqref{omegamu} is of rank three. Let
 $\varpi_\mu:\bbR^{d+1}\to\bbR^{d+1}$ be an invertible map with $\varpi_i^\mu=\omega_i^\mu$ for
$i=1,2,3$.
Let $g_\nu= |\det(\varpi^\mu)| f_\nu(\varpi^\mu\cdot)$.
Then the Fourier transform of $g_\nu$ is supported in
$R_\nu\times \bbR^{d+1-3}$ where the $R_\nu$ are $C$-extensions of
plates in $\bbR^3$ associated to the curve $(\alpha, \alpha^2/2)$.
Thus we can apply Wolff's theorem in
 three dimensions,
in the form  of Proposition  \ref{generalcones},  and obtain the estimate
$$\Big\|\sum_\nu g_\nu\Big\|_p\le \cC_{p,\eps}
\Big(\frac{r_0}{2^n r_1}\Big)^{1-\frac 4p+\eps}
\Big(\sum_\nu\big\|g_\nu\big\|_p^p\Big)^{1/p}
$$
where the constant does not depend on $\mu$. The
assertion follows by rescaling.
\end{proof}

\begin{lemmasub}\label{r1estimate}
Suppose $p\ge 2$,  $r_1\ge 2^{-k/2}$.
%, that $a_\nu
%\in \fS_k(r_1,s_\nu)$ and that the $s_\nu$ are $r_1$-separated.
Then
\begin{equation}\label{eqtrivialamunestimate1}
\Big(\sum_{\nu}\big\|\cT [a_{0,\nu}^\mu, f]\big\|_p^p\Big)^{1/p}
\le C r_1\|f\|_p,
\end{equation}
(with the usual $\sup$ modification for $p=\infty$).
\end{lemmasub}

\begin{proof}

We prove \eqref{eqtrivialamunestimate1} by interpolation and
it suffices to check the cases $p=2$ and $p=\infty$.
For $p=\infty$ the assertion follows from Lemma
\ref {varyings} and Lemma
\ref{lemmapointwiseKsestimate} if we observe that an additional $s$
integration is extended over an interval of length $\cO_1(r_1)$.
For  $p=2$ we use van der Corput's lemma in the $s$ variable
with two derivatives to take advantage of our nondegeneracy hypothesis.
Fix $\tau$, $\xi$  and observe that the amplitude of the oscillatory
integral (as a function of $s$) is bounded and has an integrable
derivative,
with uniform bounds.
Thus
$$\big|\fm[a^\mu_{0,\nu}](\xi,\tau)\big|\le C 2^{-k/2}$$
Observe that $|U_\mu|\le 2^{k}r_1^2$
on the support of
$a^\mu_{0,\nu}$;
 moreover, the supports of the $a^\mu_{0,\nu}$ are  essentially disjoint,
by (vii) of Lemma \ref{supportproperties}.
We obtain by Plancherel's theorem
that
\begin{align*}
\sum_{\nu}\big\|\cT [a^\mu_{0,\nu}, f]\big\|_2^2&=c
\sum_{\nu}\int_\xi\int_{\tau:|U_\mu(\xi,\tau)|\le 2^{k}r_1^2}
\big|\fm[a^\mu_{0,\nu}](\xi,\tau)\big|^2 d\tau | \widehat f(\xi)|^2 d\xi
\\&\lc 2^{-k}2^{k}r_1^2 \|f\|_2^2
\end{align*}
which is the desired bound for $p=2$.
\end{proof}

\begin{lemmasub}\label{trivialamunestimate}
For $n\ge 1$, $2\le p\le\infty$,
\begin{equation}\label{eqtrivialamunestimate2}
\Big(\sum_{\mu,\nu}\big\|\cT [a^\mu_{n\nu}, f]\big\|_p^p\Big)^{1/p}\le
C 2^{-k-n}r_1^{-1}
\|f\|_p.
\end{equation}
\end{lemmasub}

\begin{proof}
We argue similarly as in the proof of Lemma
\ref{r1estimate}
but begin by integrating   by parts
 with respect to
$t$ to get
$$\fm[a^\mu_{n,\nu}](\xi,\tau)=\iint e^{-it(\tau+\inn{\ga(s)}{\xi})}
\frac{\partial_t a^{\mu}_{n,\nu}(s,t,\xi,\tau)}
{i(\tau+\inn{\ga(s)}{\xi})} ds dt.
$$
Now expand
$\inn{\gamma(s)}{\xi}$ about $\scr(\xi)$ and by
 \eqref{approximationtwo}
\begin{align*}\tau+\inn{\gamma(s)}{\xi}&=
U_\mu(\xi,\tau)+
%\tau+\inn{\gamma(\scr(\xi))}{\xi}+
\cO_1(2^k(s-\scr(\xi)^2)
+\cO_1(2^{k+10} Mr_0^3)
\end{align*}
 in the  support of
$a^\mu_{n,\nu}$. Since $n\ge 1$ one also has
 $|s-\scr(\xi)|^2\le 2^{-k-8}|U_\mu|$ and hence
$$|\tau+\inn{\gamma(s)}{\xi}| \ge \frac 12 |U_\mu(\xi,\tau)|
\ge 2^{k+2n-4}r_1^2.$$
%  Moreover $\inn{\gamma'(s)}{\xi}=O((s-\scr(\xi)))|\xi|=
%\cO_1(2^{k+n+2 }r_1).$
Consequently, the multiplier $(\xi,\tau)\mapsto \partial_t a^{\mu}_{n,\nu}(s,t,\xi,\tau)
(\tau+\inn{\ga(s)}{\xi})^{-1}$ can be written as
$C2^{-(k+2n)}r_1^{-2}$ times a multiplier
in $\cS_k(r,s)$ and thus  Lemma
\ref{lemmapointwiseKsestimate} applies. Since we perform
 an $s$-integration over an interval
of length $O(2^nr_1)$ we get the asserted  $\ell^\infty(L^\infty)$ bound.

For the $\ell^2(L^2)$ estimate we  apply van der Corput's Lemma
with second derivatives and check using the support properties
of $a^\mu_{n,\nu}$
that the $L^\infty$ norm of the amplitude and the $L^1$ norm
(in $s$) of its
derivative is bounded by $C2^{-k+2n}r_1^{-2}$.
Thus we now obtain
\begin{align*}
\sum_{\mu,\nu}\big\|\cT [a^\mu_{n,\nu}, f]\big\|_2^2&\lc
\sum_{\mu,\nu}\int_\xi\int_{|U_\mu(\xi,\tau)|\le 2^{k+2n}r_1^2}
|2^{-k/2} (2^{k+2n}r_1^2)^{-1}|^2 d\tau |\widehat f(\xi)|^2 d\xi
\\&\lc  \big(2^{-k-n} r_1^{-1} \|f\|_2\big)^2.
\end{align*}
\end{proof}

\begin{lemmasub}\label{bmunestimate}
For $n\ge 1$, $2\le p\le\infty$,
\begin{equation}
\Big(\sum_{\mu,\nu}\big\|\cT[b^\mu_{n,\nu}, f]\big\|_p^p\Big)^{1/p}
\le C 2^{-k-n}r_1^{-1} \|f\|_p.
\end{equation}
\end{lemmasub}

\begin{proof}
We argue similarly as in Lemma
\ref{trivialamunestimate}. Now we  integrate by parts in $s$ to see that
$$
\fm[b^\mu_{n,\nu}](\xi,\tau)=\iint e^{-it(\tau+\inn{\ga(s)}{\xi})}
c^\mu_{n,\nu}(s,t,\xi,\tau)
ds dt
$$
where
$$c^\mu_{n,\nu}(s,t,\xi,\tau)=
\frac{\partial_s  b^\mu_{n,\nu}(s,t,\xi,\tau)}{it\inn{\gamma'(s)}{\xi}}
+
\frac{b^\mu_{n,\nu}(s,t,\xi,\tau)
\inn{\ga''(s)}{\xi}} {it(\inn{\gamma'(s)}{\xi})^2}.
$$
Now
$$|\inn{\gamma'(s)}{\xi}|\approx 2^k |s-\scr(\xi)|
\approx 2^k 2^n r_1$$
on the support of $b^\mu_{n,\nu}$.

The multiplier $(\xi,\tau)\mapsto c^{\mu}_{n,\nu}(s,t,\xi,\tau)$ is
$C2^{-(k+2n)}r_1^{-2}$ times a multiplier
in $\cS_k(r,s)$. Thus Lemma
\ref{lemmapointwiseKsestimate} applies and the $\ell^\infty(L^\infty)$ estimate follows in the same way as in Lemma
\ref{trivialamunestimate}.

For the $L^2$ estimate we may integrate by parts in $t$ and obtain
$$
|\fm[b^\mu_{n,\nu}](\xi,\tau)|\lc\iint
\min\{1,|\tau+\inn{\gamma(s)}{\xi}|^{-1}\}
|c^\mu_{n,\nu}(s,t,\xi,\tau)|dt ds,
$$
and consequently
$\sum_{\mu,\nu}\|\cT [b^\mu_{n,\nu}, f]\|_2^2$ is dominated by
\begin{align*}
&\sum_{\mu,\nu}\iint_{supp (b^\mu_{n,\nu})} \Big[
\int\limits_{|s-\scr(\xi)|\atop {\le 2^nr_1}}
\frac{2^{-k-2n} r_1^{-2}}{1+|\tau+\inn{\gamma(s)}{\xi} |}
ds\Big]^2\, d\tau
 |\widehat f(\xi)|^2 d\xi
\\
&\le 2^{n+1} r_1 \sup_\xi\Big[
\int\limits_{|s-\scr(\xi)|\atop {\le 2^nr_1}}\int_\tau
\frac{2^{-2k-4n} r_1^{-4}}{(1+|\tau+\inn{\gamma(s)}{\xi} |)^2}
d\tau ds\Big]^2
\sum_{\mu,\nu}\iint_{supp (b^\mu_{n,\nu})} \int|\widehat f(\eta)|^2 d\eta
\end{align*}
which is bounded by $( 2^{-k-n} r_1^{-1} \|f\|_2)^2$.
An interpolation yields the claimed inequality for $2\le p\le\infty$.
\end{proof}

\medskip

\noi{\bf Conclusion of the proof of Proposition \ref{maininductionstep}.}
It is immediate from the decomposition
\eqref{amudecomposition}, and
Lemma \ref{wolffestimateinRd} that for all $\eps\in (0,1)$, $p>p_W$,
\begin{multline*}
\Big(\sum_\mu\big\| \cT[a^\mu,f]\big\|_p^p\Big)^{1/p}
 \le \cC_{p, \eps} (r_0/r_1)^{1-\frac 4p +\eps}
\Big[
\Big(\sum_\nu\big\| \cT[a^\mu_{0,\nu},f]\big\|_p^p\Big)^{1/p} \\ +
\Big(\sum_{n\ge 1\atop{2^n r_1\le r_0}}
\sum_\nu\big\| \cT[a^\mu_{n,\nu},f]\big\|_p^p +
\big\| \cT[b^\mu_{n,\nu},f]\big\|_p^p
\Big)^{1/p}\Big]
\end{multline*}
and we can apply
 Lemma
\ref{trivialamunestimate} and Lemma \ref {bmunestimate} to the terms
 involving $n\ge 1$. Thus the left hand side of the inequality is
dominated by
$$
\cC_{p,\eps} (r_0/r_1)^{1-\frac 4p +\eps} \Big[
\Big(\sum_\nu\big\| \cT[a^\mu_{0,\nu},f]\big\|_p^p\Big)^{1/p} +
2^{-k}r_1^{-1}\|f\|_p\Big].
$$
There is a constant $C$ so that the symbols
$C^{-1}a^{\mu}_{0,\nu}$
belong to $\fS_k(r_1,s_{0\nu})$. Moreover, given fixed $\nu$ the function
 $a^\mu_{0,\nu}$ is not identically $0$ for at most
five $\mu$ and the $s_{0\nu}$ are $r_1$-separated. By
 a pidgeonhole  argument
we deduce the assertion of the proposition.
\qed

\medskip

\noi{\bf Conclusion of the proof of  Theorem \ref{micls}.}
Let  $\eps_1= \eps_0^2/(dM)$.
For fixed  $s_0$ and we shall  prove an estimate
for the symbol $\widetilde b(s,t,\xi)=
b(s,t,\xi) \chi(2^{k\eps_1}(s-s_0))$, namely
\begin{equation}\label{fBAestimate}
\|\fA_\gamma[\widetilde b,f]\|_{L^p(\bbR^{d+1})}\le C_p(\eps_0)
2^{-k(\frac 2p-\frac{\eps_0}2)}\|f\|_p, \quad p>p_W.
\end{equation}
The assertion of te theorem follows from
\eqref{fBAestimate}
 as
$b$ is a sum of $O(2^{k\eps_1})$ such symbols.

We write by using the Fourier inversion formula in $\bbR^{d+1}$
$$\fA_\gamma[\widetilde b,f]= \sum_{\ell=0}^\infty
\cT[\widetilde a_\ell,f]$$
where $$\widetilde a_0(s,t,\xi,\tau)=
\widetilde b(s,t,\xi) \eta_0( 2^{2k\eps_1}
(2^{-k}|\tau+\inn{\gamma(\scr(\xi))}{\xi}|+(s-\scr(\xi))^2))
$$ and, for $\ell\ge1$,
$$\widetilde a_\ell(s,t,\xi,\tau)=
\widetilde b(s,t,\xi) \eta_1( 2^{2k\eps_1-2\ell+2}
(2^{-k}|\tau+\inn{\gamma(\scr(\xi))}{\xi}|+(s-\scr(\xi))^2)).
$$

We first show the  main estimate which  is
\begin{equation}\label{mainboundinlsm}
\|\cT[\widetilde a_0 ,f]\|_{L^p(\bbR^{d+1})}
\le C_{p}(\eps_0) 2^{-k (2/p -\eps_0/2)} \|f\|_{L^p(\bbR^d)}, \quad p>p_W.
\end{equation}

Now for a constant $C$ we have $C^{-1} \widetilde a_0\in
\fS_k( 2^{-k\eps_1}, s_0)$.
We  apply Proposition
\ref{maininductionstep} iteratively
choosing $r_0$, $r_1$ to be
\begin{equation*}
r_0(n)= (2^{-k\eps_1}M)^{(3/2)^n},\quad
r_1(n)= (2^{-k\eps_1}100M)^{(3/2)^{n+1}},
\end{equation*}
for  $n=0,\dots, N$, where $N=N(\eps_1)$ is the largest integer for which
$r_1(n)\ge 2^{-k(\frac 12-\eps_1)}$. Thus certainly
%\begin{equation}\label{roneNbound}
$r_1(N)
\le 2^{-\frac k2 +2k\eps_1}$
%\end{equation}
and
%\begin{equation}\label{boundonN}
$N=N(\eps_1)\le C/\eps_1\le C'/\eps_0^2.$
%\end{equation}

By Proposition
\ref{maininductionstep} we obtain for all $p>p_W$, $\eps>0$ that
\begin{multline}\label{applicationofprop}
\big\|\cT[\widetilde a_0,f]\big\|_p
\le (\cC_{p,\eps})^{N}
 r_1(N)^{-(1-4/p+N\eps)}\Big(\sum_\nu\|\cT[a_\nu,f]\|_p^p\Big)^{1/p}
\\
+2^{-k} \sum_{n=0}^{N} (\cC_{p,\eps})^{n}r_1(n)^{4/p-2-n\eps}
 \|f\|_p
\end{multline}
where $a_\nu\in \fS_k(s_\nu,r_1(N))$ and the $s_\nu$ are $r_1(N)$-separated
points.
Note that since $r_1(N)\ge 2^{-k/2}$
\begin{equation} \label{geometricsum}
2^{-k} \sum_{n=0}^{N} r_1(n)^{4/p-2-n\eps}\le
C2^{-k(\frac 2p-N\eps)} (2^kr_1^2)^{-(1-\frac 2p)};\end{equation}
moreover  by Lemma \ref{r1estimate}
\begin{equation}\label{applicationoflemma}
\Big(\sum_\nu\big\|\cT[a_\nu,f]\big\|_p^p\Big)^{1/p}\lc r_1(N) \|f\|_p.
\end{equation}

We choose $\eps=(\eps_1/(1000 C))^2$ in
\eqref{applicationofprop} which makes $N\eps\ll \eps_0/4$ and still,
by our
previous  choice of $\eps_1$, the resulting constant
$(\cC_{p,\eps})^{N}$ depends only on $\eps_0$ and $p$.
We combine the resulting bound with
 \eqref{geometricsum} and
\eqref{applicationoflemma} and
the main estimate \eqref {mainboundinlsm} follows.

To finish the proof we have to dispose of the terms
$\cT[a_\ell, \cdot]$ for $\ell\ge 1$; these are  error terms which
can be handled by standard arguments.
We split (in analogy to a previous decomposition)
$\widetilde a_\ell=\widetilde a_{\ell,1}+\widetilde a_{\ell,2}$
where
\begin{align*}
\widetilde a_{\ell,1}(s,t,\xi,\tau)&=\widetilde a_{\ell} (s,t,\xi,\tau)
\eta_0\big(M\frac{(s-\scr(\xi))^2}{
2^{-k}|\tau+\inn{\gamma(\scr(\xi))}{\xi}|}\big),
\\
\widetilde a_{\ell,2}(s,t,\xi,\tau)&=\widetilde a_{\ell}
(s,t,\xi,\tau)
\big(1-
\eta_0\big(M\frac{(s-\scr(\xi))^2}{
2^{-k}|\tau+\inn{\gamma(\scr(\xi))}{\xi}|}\big)\big);
\end{align*}
note that $\widetilde a_{\ell,2}=0$ if $\ell \gg 2^{k\eps_1}$.
We use an integration by parts in $t$ for the integral defining
$\fm[\widetilde a_{\ell,1}]$ and
an integration by parts in $s$ for the integral defining
$\fm[\widetilde a_{\ell,2}]$.

Now for $i=1,2$
$$
\cT[\widetilde a_{\ell,i}, f](x,t')= \int \int
K_s [m_{s, \ell,i}] (x-y,t') ds\, f(y) dy
$$
where we use the notation \eqref{Ksm} with
\begin{align*}
m_{s,\ell,1}(t,\xi,\tau)&=
\frac{\widetilde a_{\ell,1}(s,t,\xi,\tau)}
{i(\tau+\inn{\gamma(s)}{\xi})},
\\
m_{s,\ell,2}(t,\xi,\tau)&=
\Big[
\frac{\partial_s  \widetilde a_{\ell,2}(s,t,\xi,\tau)}{it\inn{\gamma'(s)}{\xi}}
+
\frac{ \widetilde a_{\ell,2}(s,t,\xi,\tau)
\inn{\ga''(s)}{\xi}} {it(\inn{\gamma'(s)}{\xi})^2}\Big].
\end{align*}
We argue as in the proof of  Lemma \ref{lemmapointwiseKsestimate}
and by a straightforward integration by parts we  obtain the bounds
\begin{multline*}
|K_s[m_{s,\ell,i}](x,t')|
\\
% C \int
%(1+2^k r^2|t-t'|)^{-2}\quad\times
 C 2^{-k(1-2\eps_1)-\ell}\int
\frac{2^{kd} 2^{\ell}}
{(1+2^{k(1-2\eps_1)+\ell} |t-t'|+
2^{k(1-2\eps_1) }|x-t'\gamma(s)|)^{d+2}}\,
 dt
\end{multline*}
for $i=1,2$; here we use for the second kernel  that $m_{s,\ell, 2}=0$ for
$\ell \ge 2^{2k\eps_1}.$

This estimate implies (after an integration in $s$)
 that the terms involving $\widetilde a_\ell$ for $\ell>0$
are error terms and we get the estimates
$$\|\cT[\widetilde a_\ell ,f]\|_{L^p(\bbR^{d+1})}
\lc 2^{-k(1-2\eps_1 (d+1))} \|f\|_{L^p(\bbR^d)}
$$
for $1\le p\le \infty$ and of course the constant here is much smaller
 than $2^{-2k/p}$ for $p>4$ and in particular for $p>p_W$. This finishes the proof of Theorem
\ref{micls}.
\qed

\subsection{Proof of Theorem \ref{genmicls}}
We may assume that $B>100(d+\delta^{-1}) $. In addition by a reparametrization
we  may also assume that
$\Gamma$ is parametrized by arclength $s$ (consequently we may have to replace $B$ by a power of $B$).

We localize in $s$ (splitting the parameter interval in
$O(B^{102})$ pieces) and assume that the symbol  is localized to an
$s$-interval $I(s_0)$ centered at $s_0$, and of length $\le B^{-100}$.
By further localization in $\Omega$ we split the symbol
into $O(B^{100 d})$ pieces
localized in balls of the form
$\Omega(\xi^0)=\{\xi:|\xi-\xi^0|\le 2^{k} B^{-10}\}$ where
$B^{-1}2^k\le |\xi^0|\le 2^kB$. We now assume that our symbol $a$
 is
supported in $I(s_0)\times[1,2]\times \Omega(\xi^0)$ and that $a$
satisfies differentiability conditions similar to
\eqref{aassumption}, but with the  constant
$\cC[a] $ replaced by
$C_d  \cC[a] B^{1000 d}$; moreover we assume the lower bound
\begin{equation}\label{Gammasecondlowerboundmod}
\big|
\inn{\Gamma'(s)}{\xi}\big |+
\big|\inn{\Gamma''(s)}{\xi}\big |
\ge 2B^{-2}|\xi| \quad \text{ if  } (s,t,\xi)\in \supp(a).
\end{equation}

We set $\theta=\xi^0/|\xi^0|$ and  distinguish two cases.
In the first case we assume that
$|\inn{\Gamma'(s_0)}{\theta^0}| \ge B^{-100}$; then by the support
properties after localization
% we have
$\inn{\Gamma'(s)}{\xi}| \ge B^{-90}|\xi|$ on the support of $a$.
This allows us to perform an integration by parts in $s$ first,
 thus gaining a power of  $2^k$  and standard estimates yield
that in the present case the $L^p(\bbR^d)$ norm of
$\fA_\Gamma[a,f](\cdot,t)$ is bounded by $C_{B,d} 2^{-k}\|f\|_p$
for $1\le p\le \infty$, uniformly in $t\in [1,2]$. Thus in this
case we obtain a better bound than the one claimed in
\eqref{Aafinequality}.

For the second (main) case we  have the inequalities
\begin{align}  \label{Gammaprimesmallness}
|\inn{\Gamma'(s_0)}{\theta^0}| &\le B^{-100},
%\end{equation}\begin{equation}
\\
 \label{Gammadoubleprimelargeness}
|\inn{\Gamma''(s_0)}{\theta^0}| &\ge B^{-2}.
\end{align}

Now let
 $\{v_1,\dots, v_d\}$ be an orthonormal basis of $\bbR^d$ so that $v_1= \Gamma'(s_0)$ and $\text{span }\{v_1,v_2\}=\text{span}
\{\Gamma'(s_0),\theta^0\}$. Let $L$ be the linear transformation with
$L(v_i)=v_i$ for $i=1$, $3\le i\le d$ and $L(v_2)=\inn{\Gamma''(s_0)}{v_2}^{-1} v_2$.
Let $\widetilde \Gamma(s)= L\Gamma(s)$.

By
\eqref{Gammaprimesmallness} and
\eqref{Gammadoubleprimelargeness}
and a Taylor expansion
 $|\widetilde \Gamma'(s)|=1+O(B^{-10})$
and  since we assume that $\Gamma$ is parametrized by arclength
we have $\inn{\Gamma'(s)}{\Gamma''(s)}=0$. A calculation shows that
$\inn{\widetilde \Gamma''(s)}{\theta^0}=1 +O(B^{-10}).$

Notice that
\begin{equation}\label{Lscaling}\fA_\Gamma[a,f](x,t)=
\fA_{\widetilde \Gamma}[\widetilde a, f\circ L^{-1}] (Lx,t)
\end{equation}
where $\widetilde a (s,t,\eta)=a(s,t, L^t\eta)$.
After a reparametrization of $\widetilde \Gamma$ by arclength
 an application of Theorem \ref{micls} shows
that
% Then for $p>p_W$
\begin{equation*}
\big\|\fA_{\widetilde \Gamma}[\widetilde a,f]\big\|_{L^p(\bbR^{d+1})}
\le C(\eps_0,p, B) \cC[a] 2^{-k (\frac 2p-\eps_0)}
 \|f\|_{L^p(\bbR^{d})}, \quad p>p_W,
\end{equation*}
and the corresponding assertion
for $\fA_\Gamma$ follows by  \eqref{Lscaling}.\qed

\section{Local Smoothing for Curves in $\bbR^3$}\label{LSM2}

We now return to the situation in $\bbR^3$ and consider
\emph{curves with nonvanishing curvature and torsion.} We
shall use notation as in \S\ref{averages} and prove an estimate
for the $t$-dilates of the operators $\cA^{k,l,\nu}$ defined in
\eqref{defAk}. The following lemma is proved by rescaling
and the results of the previous sections.
Define
\begin{equation}
\label{Aklnt}
\widehat {\cA^{k,l,\nu}_t f}(\xi)= m_k[a_{k,l,\nu}(t\cdot)]\widehat f(\xi)\end{equation}
and let $\chi$ be a smooth function supported in $(1/2,2)$.
%Similarly we can define the $t$-dilates  $\widetilde \cA_t^{k,\nu}$,
%$\cB_t^{k,l,\nu}$, $\widetilde \cB_t^{k,\nu}$ of the operators  in
%\eqref{defAk}, \eqref{defAb}; however for these operators the fixed time
%estimates of
%\S  will suffice.

\begin{proposition}\label{rescalingtonondegeneratecase}
For $p>p_W$, $l<k/3$, $\eps>0$,
\begin{equation*}
\Big(\int \|\chi(t)\cA_t^{k,l,\nu} f\|_p^p dt\Big)^{1/p}
\le C_\eps 2^{-l(1-6/p)} 2^{-2k/p}2^{k\eps} \|f\|_p.
\end{equation*}
\end{proposition}
\begin{proof}
The symbol $a_{k,l,\nu}$ in
\eqref{Aklnt} is supported in a set where
$|\inn{\xi}{T(s_\nu)}|\approx 2^{k-2l}$,
$|\inn{\xi}{N(s_\nu)}|\lc 2^{k-l}$,
$|\inn{\xi}{B(s_\nu)}|\approx 2^{k}$.
We shall rescale the parameter $s=s_\nu +2^{-l} u$ with $u\lc 1$,
moreover we rescale in $\xi$ as follows. Let $U_\nu$ be the rotation
which maps the unit vectors $e_1,e_2,e_3$ to
$T(s_\nu), N(s_\nu), B(s_\nu)$. Let
$\Delta_l\eta=(2^l\eta_1, 2^{2l}\eta_2, 2^{3l}\eta_3)$ and let
$L_{l,\nu}=U_\nu\circ \Delta_l$. Then
% the preimage of the support of
$$(u,\eta)\to c_{k,l,\nu}(u,\eta):=a_{k,l,\nu}(s_\nu+2^{-l}u, L_{l,\nu}\eta)$$
is supported in a set where  $|\eta|\approx 2^{k-3l}$ and $|u|\lc 1$ and
there are the estimates
$$
\big|\partial_u^{(n)}\partial_\eta^{\alpha} c_{k,l,\nu}(u,\eta)\big|
\le C_{n,\alpha} 2^{-(k-3l)|\alpha|}.
$$
Moreover if we set
$$
\Gamma_{l,\nu}(u)= L_{l,\nu}^*\gamma(s_\nu+2^{-l}u)
$$
then $\Gamma_{l,\nu}$ is a $C^5$ curve
with upper bounds uniformly in $l,\nu$ and we also have
$$|\inn{\Gamma_{l,\nu}''(u)}{\eta}|\approx|\eta|\approx 2^{k-3l}$$
in the support of $c_{k,l,\nu}$ (again with the implicit constants uniform
 in $\ell,\nu$).

Changing variables we get

\begin{align} (2\pi)^d \cA_t^{k,l,\nu} f(x)&=2^{-l}
\iint e^{i t\inn{\gamma(s_\nu+2^{-l}u)}{\xi}+i\inn{x}{\xi}}
a(s_\nu+2^{-l}u, \xi) \widehat f(\xi) d\xi du
\notag
\\
&=2^{-l}
\int e^{i t\inn{\Gamma_{l,\nu}(u)}{\eta}+i\inn{L_{l,\nu}^*x}{\eta}}
c_{k,l,\nu}(\eta) \widehat f(L_{l,\nu}\eta) \,2^{6l} d\eta\, du
\notag
\\ \label{conjugationwithLlnu}
&= 2^{-l}  T_t^{k,l,\nu} [f({L_{l,\nu}^*}^{-1}\cdot)]
 (L_{l,\nu}^*x)
\end{align}
where
\begin{equation*}
T_t^{k,l,\nu}g(x) =
\iint e^{it\inn{\Gamma_{l,\nu}(u)}{\eta}} c_{k,l,\nu}(u,\eta) \widehat g(\eta) e^{i\inn x\eta} d\eta du.
\end{equation*}
Thus  we can apply Theorem \ref{genmicls}
%with the modification in
% Remark \ref{remarkonconstants},
 for the dyadic annulus of width $2^{k-3l}$,
and obtain
$$\Big(\int\big\|\chi(t) T_t^{k,l,\nu} g\big\|_p^p dt\Big)^{1/p}
\le C_\eps  2^{-2(k-3l)/p} 2^{(k-3l)\eps} \|g\|_p.
$$
We rescale using \eqref{conjugationwithLlnu} to
obtain the asserted bound.

\end{proof}

\noi{\bf Proof of Theorem \ref{localsmoothing}.}
We apply  inequality  \eqref{Step1}, rescaled by the  $t$-dilation,
  and combine it with Proposition
 \ref{rescalingtonondegeneratecase} to obtain
\begin{align*}
\Big(\int \big\|\chi(t)\sum_\nu \cA_t^{k,l,\nu} f\big\|_p^pdt\Big)^{1/p}
&\le C_\eps 2^{2l(\frac 12-\frac 2p+\eps)}
\Big(\sum_\nu \int\big\|\chi(t)
\cA^{k,l,\nu}_t f\big\|_p^p dt\Big)^{1/p}
\\&\le C_\eps' 2^{-k(\frac 2p-\eps)} 2^{-2l(\frac 1p-\eps)} \|f\|_p
\end{align*}
and thus
\begin{equation}\label{mainestimateinpfofthm14}
\Big(\int \big\|\chi(t) \cA_t^{k,l} f\big\|_p^pdt\Big)^{1/p}
\le C_\eps 2^{-2(k-l)/p +2k\eps}\|f\|_p, \quad p>p_W.
\end{equation}
This is the main estimate  and we may sum over $l<k/3$.
There are similar estimates for the operators
$\sum_\nu \widetilde \cA^{k,\nu}_t$ and $\sum_\nu\cB_t^{k,l,\nu} $
obtained if we scale by $t$ in the definitions \eqref{defAk},
\eqref{defAb}; however these follow already by integrating out the
 fixed time estimates  implied by
Proposition \ref{besov}. The conclusion is that if $m_k$ is as in
\eqref{dyadicmultipliers} then
$$\Big(\int\|\chi(t) \cF^{-1}[m_k(t\cdot )\widehat f]\|_p^p dt
\Big)^{1/p}
\le C_\eps 2^{-k (\frac {4p}3-\eps)}\|f\|_p$$
and the assertion of Theorem  \ref{localsmoothing} on boundedness in
Sobolev spaces
follows by standard arguments.\qed

\section{Maximal Functions}\label{maximalfunctions}

\noi {\bf Proof of Theorem \ref{maximal}.}
Given  Theorem \ref{localsmoothing} the proof  is straightforward for the
case of  curves with nonvanishing curvature and torsion.
Let $\cL_k$ a  Littlewood-Paley operator which localizes  to frequencies
of size
$\approx 2^k$.
Then for $p>p_W$
\begin{equation*}
%\label{locsmzeroderivative}
\Big(\int_1^2 \big\|\cA_{2^\ell t} \cL_{k+\ell} f\big\|_p^p dt\Big)^{1/p}
\le C_{\eps,p} 2^{-k(\frac 4{3p}-\eps)}  \|\cL_{k+\ell} f\|_p,
\end{equation*}
and
\begin{equation*}
%\label{locsmfirstderivative}
\Big(\int_1^2 \big\|(\partial/\partial t)
\cA_{2^\ell t} \cL_{k+\ell} f\big\|_p^p dt\Big)^{1/p}
\le C_{p,\eps} 2^{-k(\frac 4{3p}-\eps)}(1+2^k \sup_{s\in I} |\gamma(s)|)
\|\cL_{k+\ell} f\|_p,
\end{equation*}
and by standard arguments  we obtain
\begin{multline*}
%\label{maximalfunctionestimate}
\big\|\sup_{\ell\in \bbZ}\sup_{1\le t\le 2} |\cA_{2^\ell t} \cL_{k+\ell}f|\big\|_p\\
\le
C_{p,\eps} 2^{-k(\frac 4{3p}-\eps)}\big (1+2^k \sup_{s\in I}
|\gamma(s)|\big)^{1/p}
\Big( \sum_{\ell\in\bbZ} \|\cL_{k+\ell} f\|_p^p\Big)^{1/p}, \quad p>p_W.
\end{multline*}
Since $p\ge 2$ we have
$( \sum_{\ell\in\bbZ} \|\cL_{k+\ell} f\|_p^p)^{1/p} \lc\|f\|_p$.
Similar  $L^2$ estimates based on van der Corput's lemma yield
\begin{equation*}
%\label{maximalfunctionestimateL2}
\big\|\sup_{\ell\in \bbZ}\sup_{1\le t\le 2} |\cA_{2^\ell t}\cL_{k+\ell}
 f|\big\|_2
\le
C 2^{-k/3}\big (1+2^k \sup_{s\in I} |\gamma(s)|\big)^{1/2}\|f\|_2
\end{equation*}
and an interpolation shows that
\begin{equation}\label{interpolatedmaximalestimate}
\big\|\sup_{\ell\in \bbZ}\sup_{1\le t\le 2} |\cA_{2^\ell t} \cL_{k+\ell}f|
\big\|_p
\lc C_p (1+\sup_{s\in I}|\gamma(s)|)^{1/p} 2^{-k a(p)}\|f\|_p
\end{equation}
with $a(p)>0$ if $p>(p_W+2)/2$.
This proves the statement of Theorem \ref{maximal}
in the case of curves with nonvanishing curvature and torsion.

In the finite type case we use rescaling as in \S\ref{averages}.
We
 may after using a partition of unity assume that
\eqref{expansionofgamma} holds, with $n_1<n_2<n_3$, and $n_3\le n$.
Then, with $A_{j,t}$ as in  \eqref{definitionofAtj}
we need to show that
\begin{equation}\label{claimedscaledmaximalestimate}
\big\|\sup_\ell\sup_{1\le t\le 2}|A_{j,2^\ell t}f| \big\|_p
\lc 2^{j(\frac{n_3}{p} -1)} \|f\|_p, \quad p>(p_W+2)/2.
\end{equation}
We may apply \eqref{interpolatedmaximalestimate}
to the normalized curves $\delta_j\gamma(s_0)+\Gamma_j$
(where $\Gamma_j$ is as in \eqref{finitetypedilation}), and observe that
$$\sup_u|\delta_j\gamma(s_0)+\Gamma_j(u)|=O(2^{n_3}).$$
Thus using also
\eqref{finitetypesubstitution}
and setting %we get with
$f_{-j}=f\circ \delta_{-j}$, $N=n_1+n_2+n_3$, we see that
$\|\sup_\ell\sup_{1\le t\le 2}|A_{j,2^\ell t}f| \|_p$  is controlled by
\begin{align*}
%&\big\|\sup_\ell\sup_{1\le t\le 2}|A_{j,2^\ell t}f| \big\|_p
%\\
%&\lc 2^{-j}\Big\|
& 2^{-j}\Big\|
\sup_\ell\sup_{1\le t\le 2}\big|
\int f_{-j}(\delta_j \cdot -2^\ell t\delta_j\gamma(s_0)-2^\ell t\Gamma_j(u))
\chi_j(u) du\big|\, \Big\|_p
\\
& \lc 2^{j(\frac{n_3}{p} -1)} 2^{-N/p}\|f_{-j}\|_p
\lc 2^{j(\frac{n_3}{p} -1)} \|f\|_p,
\end{align*}
and  obtain \eqref{claimedscaledmaximalestimate}. We need to sum in $j$
in \eqref{claimedscaledmaximalestimate}
 which is possible since also $p> n_3$.
\qed

\medskip

\noi{\bf A two-parameter maximal function.} Our results on local smoothing
can also be used
to  prove  bounds for certain  two-parameter
maximal functions. Consider the  two-parameter family of helices
\[ H(a,b):= \left\{\gamma_{a,b}=(a \cos (2\pi s),\, a \sin(2\pi  s),\, bs) \,:\, 0 \leq s \leq 1
\right\}, \, 1<a,b<2. \]
Then we obtain a lower bound on the Hausdorff dimension of some ``Kakeya-type'' sets.

\begin{proposition}  \label{Kakeya}
Let $F$ be a set which for every ${x} \in \mathbb R^3$ contains a
helix $x + H(a,b)$ for some $(a,b)$, $1<a,b<2$.
Then the Hausdorff dimension of $F$
is at least 8/3.
\end{proposition}

By arguments in \cite{Bourgain1991} one sees
 that Proposition \ref{Kakeya}
is a consequence of an estimate for a local maximal operator, namely
\begin{equation}\label{2parametermaximal}
\Big\|\sup_{1\le a,b\le 2}\big|\int f(x-\gamma_{a,b}(s)) \chi(s) ds\big|
\Big\|_p \le C_\alpha \|f\|_{L^p_\alpha}, \quad p>p_W,
\alpha>(3p)^{-1}.
\end{equation}

\begin{proof}
We only sketch the argument since it  follows the same lines as the
one  in  the proof of
Theorem \ref{maximal}, however it uses as an
 additional ingredient the relation between $\partial_a\gamma_{a,b}$ and
$\gamma''$.

Let $\widetilde \cA^k_{a,b}$,  $\cA^{k,l}_{a,b}$, $\cB^{k,l}_{a,b}$
be the operators with symbols
$m_k(\widetilde a_k)$,
$m_k(a_{k,l})$ , $m_k(b_{k,l})$ as in
\eqref{definitionofwidetildeak},
\eqref{definitionofaklbkl}, for the curve $\gamma_{a,b}$.
\eqref{2parametermaximal} follows from
\begin{equation}\label{aklmaxest}
\big\|\sup_{1\le a,b\le 2}|\cA_{a,b}^{k,l} f|
\big\|_p \le C_\eps 2^{l/p} 2^{k\eps} \|f\|_{p}, \quad p>p_W,\quad l<k/3,
\end{equation}
 and  related  statements for
$\widetilde \cA^k_{a,b}$, $\cB^{k,l}_{a,b}$.
By  standard arguments the  proof of
\eqref{aklmaxest} can be  reduced to
\begin{multline}\label{differentiations}
\Big(\iint_{1\le a,b\le 2}\Big\|
\frac{\partial^{j_1+j_2}}{(\partial a)^{j_1}(\partial b)^{j_2}}
\cA^{k,l}_{a,b} f\Big\|_p^p da db\Big)^{1/p}
\\\le C_\eps 2^{(k-l)j_1+kj_2} 2^{-2(k-l)/p+k\eps} \|f\|_p,
\end{multline}
$l<k/3$.
When $j_1=0$, $j_2\in \{0,1\}$ inequality
\eqref{differentiations} follows  from
\eqref{mainestimateinpfofthm14}.
For the $a$-differentiation ($j_1=1$) inequality
\eqref{differentiations} asserts a blowup of merely $2^{k-l}$.
This happens because
the $a$-differentiation of the phase yields an additional factor of
$$\partial_a \inn{\gamma_{a,b}(s)}{\xi}=
\xi_1\cos (2\pi s)+\xi_2 \sin(2\pi  s)= -(4\pi^2a)^{-1}
\inn{\gamma_{a,b}''(s)}{\xi},
$$
for the symbol,
and $\inn{\gamma_{a,b}''(s)}{\xi}$ is of size $\approx 2^{k-l}$
on the support of $m_k(a_{k,l})$. It is here where we
use
the improvements stated in
part (iv) of  Lemma \ref{platesandestimates} and part (iii) of
Lemma \ref{L1L2bounds}.
\end{proof}

%\section{$L^p$-regularity of a restricted X-ray transform}
%\footnote{INCLUDE OR NOT?}

\end{document}